\input amstex
\documentstyle{amsppt}
\magnification=1200 \hsize=13.8cm \catcode`\@=11
\def\NoLogo{\let\logo@\empty}
\catcode`\@=\active \NoLogo

\def\tN{\widetilde N}

\def\tM{\widetilde M}

\def\heatau{\lf(\frac{\p}{\p \tau}-\Delta\ri)}
\def\taub{\bar{\tau}}
\def\sigmab{\bar{\sigma}}
\def\texp{\widetilde \exp}
\def\heat{\lf(\frac{\p}{\p t}-\Delta\ri)}
\def \b {\beta}

\def\Ric{\text{Ric}}
\def\lf{\left}
\def\ri{\right}
\def\bbar{\bar \beta}
\def\a{\alpha}

\def\g{\gamma}
\def\e{\epsilon}
\def\p{\partial}
\def\delbar{{\bar\delta}}

\def\dbar{\bar\partial}

\def\C{\Bbb C}
\def\R{\Bbb R}
\def\P{\Bbb P}

\def\vp{\varphi}

\def\dbar{\bar\partial}

\def\bb{{\bar\beta}}
\def\abb{{\alpha\bar\beta}}
\def\gbd{{\gamma\bar\delta}}

\def \D {\Delta}
\def\aint{\frac{\ \ }{\ \ }{\hskip -0.4cm}\int}
\documentstyle{amsppt}
\magnification=1200 \hsize=13.8cm \vsize=19 cm

\leftheadtext{ Lei  Ni} \rightheadtext{Li-Yau-Hamilton inequality}
\topmatter
\title{A new matrix Li-Yau-Hamilton estimate for K\"ahler-Ricci flow}\endtitle

\author{ Lei Ni\footnotemark }\endauthor
\footnotetext"$^{1}$"{Research partially supported by NSF grants
and an Alfred P. Sloan Fellowship, USA.}

\address
Department of Mathematics, University of California, San Diego, La
Jolla, CA 92093
\endaddress
\email{lni\@math.ucsd.edu}
\endemail

\affil { University of California, San Diego}
\endaffil

\date  September 2004 (Revised, March 2005) \enddate

\abstract In this paper we prove a new matrix Li-Yau-Hamilton
estimate for K\"ahler-Ricci flow. The form of this  new
Li-Yau-Hamilton estimate is obtained by the interpolation
consideration originated in \cite{Ch1}. This new inequality is
shown to be connected with Perelman's entropy formula through a
family of differential equalities. In the rest of the paper, We
show several applications of this new estimate and its linear
version proved earlier in \cite{CN}. These include a sharp heat
kernel comparison theorem, generalizing
 the earlier result of Li and Tian, a manifold version of Stoll's theorem
 on the characterization of `algebraic divisors', and a
 localized monotonicity formula for analytic subvarieties.

Motivated by  the connection between the heat kernel estimate and
the reduced volume monotonicity of Perelman,  we  prove a sharp
lower bound heat kernel estimate for the
 time-dependent heat equation, which is, in a certain sense,  dual to Perelman's
 monotonicity of the `reduced volume'. As an application of this new monotonicity formula,
 we show that the blow-down
 limit of a certain type III immortal solution is a gradient  expanding soliton.
 In the last section we also illustrate the
 connection between  the new Li-Yau-Hamilton estimate and  the
earlier  Hessian comparison theorem on  the `reduced distance',
proved in \cite{FIN}.
\endabstract

\endtopmatter

\document

\subheading{\S0 Introduction}\vskip .2cm

In \cite{LY}, Peter Li and S.-T. Yau developed the fundamental
gradient estimates for positive solution $u(x,t)$ to the heat
equation $\heat u(x,t)=0$. On a complete Riemannian manifold $M$
with nonnegative Ricci curvature, they also derived a sharp form
of the classical Harnack inequality (cf. \cite{Mo}) out of their
gradient estimates. Later in \cite{H2}, Richard Hamilton extended
the estimate of Li-Yau to the full matrix version on the Hession
of $\log u$, under the stronger assumption that $M$ is Ricci
parallel and of nonnegative sectional curvature. More recently, in
\cite{CN}, H.-D. Cao and the  author observed that if $M$ is  a
K\"ahler manifold with nonnegative bisectional curvature, one can
obtain the Hamilton's matrix version estimate on the complex
Hessian of $\log u$ without the assumption of Ricci being
parallel. Following \cite{NT1}, We called our estimate in
\cite{CN} a Li-Yau-Hamilton estimate (or inequality, which is also
referred as {\it differential Harnack inequality} in some
literatures). For Ricci flow (K\"ahler-Ricci) flow there also
exists the fundamental work of Hamilton \cite{H1} (H.-D. Cao
\cite{C1}) on the matrix Li-Yau-Hamilton inequality for curvature
tensors. See also \cite{CC1-2}, \cite{CH}, \cite{CK1}, \cite{NT1},
etc, for the later developments, \cite{A} and references therein
for the Li-Yau-Hamilton estimates of other geometric flows. The
relation with the monotonicity formulae was discussed in \cite{N2,
N4}.

In this paper we shall prove a nonlinear version of the
Li-Yau-Hamilton estimate of \cite{CN}, for time dependent K\"ahler
metrics evolving by K\"ahler-Ricci flow. The new matrix inequality
asserts that if  $(M, g(t))$ is a complete solution to
K\"ahler-Ricci flow $\frac{\p}{\p t}g_{\abb}(x,t)=-R_{\abb}(x,t)$
with (bounded, in case $M$ is not compact)  nonnegative
bisectional curvature, and if   $u$ is a positive solution to the
{\it forward conjugate heat equation}: $\left(\frac{\p}{\p t}-\D
-{\Cal R}\right)u(x,t)=0$, where ${\Cal R}$ is the scalar
curvature, then
$$
u_{\abb}+\frac{u}{t}g_{\abb}+ uR_{\abb} +u_{\a}V_{\bb} +u_{\bb}
V_{\a}+uV_{\a}V_{\bb}\ge 0 \tag 0.1
$$
for any $(1,0)$ vector field $V$. The form of this new matrix
Li-Yau-Hamilton type estimate is a natural one after we found it
through an interpolation consideration which was originated in
\cite{Ch1} by Ben Chow.  We shall illustrate this interpolation
consideration further in the following paragraph.

In \cite{CH}, Ben Chow and Richard Hamilton proved a linear trace
Li-Yau-Hamilton inequality for the symmetric  positive definite
$2$-tensors evolved by the (time-dependent) Lichnerowicz heat
equation, coupled with  the Ricci flow, whose complete solution
metrics have bounded non-negative curvature operator. This result
in particular generalizes the trace form of Hamilton's fundamental
matrix Li-Yau-Hamilton estimate on curvature tensors, for
solutions to Ricci flow in \cite{H1}. Later, in \cite{Ch1}, Chow
discovered a very interesting interpolation phenomenon. Namely he
shows a family of Li-Yau-Hamilton estimates in the case of Riemann
surfaces with positive curvature such that this family  connects
Li-Yau's estimate for the positive solutions to a heat equation
with Chow-Hamilton's linear trace estimate on the solutions to the
Lichnerowicz heat equation, coupled with  the Ricci flow,  in the
case of Riemann surfaces. Seeking the
 analogue of such interpolation in higher dimensions turns out to be
fruitful. Even though the interpolation itself has not been found
directly useful in geometric problems, it does play a crucial role
in discovering new (useful) estimates. For example in \cite{N4},
such consideration led the author to discover a Li-Yau-Hamilton
estimate for the Hermitian-Einstein flow (what proved there is
more general).  This estimate of \cite{N4} is (crucial) one of the
new ingredients in obtaining the (sharp) estimates on the
dimension of the spaces of holomorphic functions (sections of
certain line bundles) of the polynomial growth. (For more details,
please see \cite{N4}, as well as  \cite{CFYZ}.) The
Li-Yau-Hamilton inequality proved in \cite{N4} can be interpolated
with the earlier one proved by L.-F. Tam and the author in
\cite{NT1}. This interpolation consideration further suggests that
there should be a one-one correspondence, for the Li-Yau-Hamilton
estimates, between the linear case and the (nonlinear) case with
Ricci flow. Seeking the linear correspondence of the  linear trace
Li-Yau-Hamilton inequality for K\"ahler-Ricci flow, proved in
\cite{NT1},  led the correct formulation of the Li-Yau-Hamilton's
inequality in \cite{N4} for the linear case. On the other hand,
looking for the nonlinear version of matrix Li-Yau-Hamilton
estimate for the linear equation  in \cite{CN} leads us to
formulate the correct form of the (nonlinear version) matrix
Li-Yau-Hamilton inequality in this paper for the K\"ahler-Ricci
flow. The proof of this result is applying the tensor maximum
principle of Hamilton, which can not be completed without the
generous help from Professor Ben Chow on a certain crucial step.
We want to record out gratitude to him here.

As in the most other cases, the new estimate proved in this paper
is sharp since it holds equality if and only if on  expanding
K\"ahler-Ricci solitons. Its proof also makes use of  the earlier
fundamental Li-Yau-Hamilton estimate of H.-D. Cao  \cite {C1} (see
also \cite{C2}). The Riemannian version of (0.1) suggests a new
matrix differential inequality on curvature tensors, which is
different from Hamilton's one in \cite{H1}.  Please see section 5
for details. This new expression also vanishes identically on
expanding solitons. Unfortunately, we have not been able to verify
this new matrix estimates at this moment. (Please see Remark 5.2
for details.)

A little surprisingly, this new Li-Yau-Hamilton estimate for the
K\"ahler-Ricci flow can be shown to be  related to Perelman's
celebrated entropy formula for the Ricci flow \cite{P}, at least
for K\"ahler case. Again this is done through Chow's interpolation
consideration. Namely, one in fact can obtain a family of {\it
pre-Li-Yau-Hamilton equalities} (a notion suggested to us by Tom
Ilmanen) such that at one end, the trace of this {\it
pre-Li-Yau-Hamilton equality}  gives the Perelman's entropy
formula after  integration on the manifold, and at the other end
one obtains the Li-Yau-Hamilton inequality of this paper by
applying the tensor maximum principle of Hamilton. At the midpoint
one can obtain both the entropy formula for the linear heat
equation proved in \cite{N3} and the Li-Yau's estimate in
\cite{LY} for the linear heat equation. (The {\it
pre-Li-Yau-Hamilton equalities} was proved earlier by Ben Chow in
\cite{Ch2} (see also \cite{CLN}) for the {\it backward Ricci flow}
on Riemannian manifolds soon after the proof of  (0.1).) This
connection between (0.1) and Perelman's entropy formula suggests
that the new Li-Yau-Hamilton estimate proved here may be of some
fundamental importance. Please see Section 3 for the detailed
exposition on this interpolation between Perelman's entropy
formula and the new Li-Yau-Hamilton estimate. One should also
refer to the recent beautiful survey \cite{Ev1} and the stellar
notes \cite{Ev2} by Evans for the relation between entropy and the
Harnack estimates for the linear heat equation (see also \cite{N3}
for a different relation). \cite{Ev1-2} also contain many other
applications of entropy consideration in the study of PDE.

The rest of the paper is on applications of this new estimate, as
well as the corresponding linear one of \cite{CN}. The immediate
consequences include the  monotonicity of a quantity which is
called Nash's entropy, a Perelman type monotonicity (or Huisken
type in the linear case) of the `reduced volume of analytic
subvarieties' and a sharp form of Harnack estimate for positive
solutions to the {\it forward conjugate heat equation}. It also
can be applied to proved a sharp heat kernel comparison theorem
for any subvariety in complete K\"ahler manifolds with nonnegative
bisectional curvature. This generalizes a previous result of Peter
Li and Gang Tian \cite{LT},  in which the authors proved the sharp
comparison on heat kernels of algebraic manifolds, equipped  with
the induced Fubini-Study metric (also called Bergmann metric in
\cite{LT}) from the ambient $\P^m$. More precisely, we proved the
following result.

\medskip

{\it Let $M$ be a complete K\"ahler manifold with nonnegative
bisectional curvature. Let $H(x,y, t)$ be the fundamental solution
of the heat equation. Let ${\Cal V}\subset M$ be any complex
subvariety of dimension $s$. Let $K_{\Cal V}(x,y, t)$ be the
fundamental solution of heat equation on ${\Cal V}$. Then
$$
K_{\Cal V}(x,y, t)\le (\pi t)^{m-s}H(x,y, t), \text{ for any }x,\,
y \in {\Cal V}. \tag 0.2
$$
The equality implies that ${\Cal V}$ is totally geodesic. }

\medskip

A more involved application is the localized monotonicity formula
and an elliptic `monotonicity principle' for the analytic
subvarieties. The later leads to a manifold (curved or nonlinear)
version of Stoll's characterization on `algebraic divisors'. This
localization uses, substantially, the beautiful ideas from the the
study of mean curvature flow of Ecker in \cite{E1-2}. We shall
give a brief sketch on these results below.

Let ${\Cal V}$ be a  subvariety of complex dimension $s$ as above.
Denote by ${\Cal A}_{{\Cal V}, x_0}(\rho)$ the $2s$-dimensional
Hausdorff measure of ${\Cal V}\cap B_{x_0}(\rho)$. Here
$B_{x_0}(\rho)$ is the ball (inside $M$) of radius $\rho$ centered
at $x_0$. The elliptic `monotonicity principle' states the
following.

\medskip

{\it There exists $C=C(m, s)$ such that for any $\rho'\in (0,
\delta(s)\rho)$
$$
\frac{{\Cal A}_{{\Cal V}, x_0}(\rho')
(\rho')^{2(m-s)}}{V_{x_0}(\rho')}\le C(m,s) \frac{{\Cal A}_{{\Cal
V}, x_0}(\rho) \rho^{2(m-s)}}{V_{x_0}(\rho)}. \tag 0.3
$$
Here $\delta(s)=\frac{1}{\sqrt{2+4s}}$. }

\medskip

The following consequence of (0.3) is somewhat interesting.

\medskip

{\it  Let $M^m$ be a complete K\"ahler manifold with nonnegative
holomorphic bisectional curvature. Suppose that $M$ contains a
compact subvariety ${\Cal V}$ of complex dimension $s$. Then there
exists $C=C(m, s)>0$ such that  for $\delta(s) \rho\ge \rho'\gg1$,
$$
\frac{V_{x_0}(\rho)}{V_{x_0}(\rho')}\le
C\left(\frac{\rho}{\rho'}\right)^{2(m-s)}. \tag 0.4
$$
In particular,
$$
\lim_{\rho\to \infty}\frac{V_{x_0}(\rho)}{\rho^{2(m-s)}}<\infty.
$$
}

\medskip
The result sharpens the Bishop-Gromov volume comparison theorem in
the presence of compact subvarieties. When $M$ is simply-connected
the result is in fact a consequence of the splitting theorem
proved in \cite{NT2, Theorem 0.4}, via a completely different
approach. It is not clear it can be derived out of any previously
known result in the general case. For an entire analytic set
${\Cal V}$ in $\C^m$ (of dimension $s$) , one can define the
Lelong number by
$$
\nu_{\infty} ({\Cal V})=\sup_{x_0\in M}\limsup_{\rho \to
\infty}\frac{(\pi\rho^2)^{(m-s)} {\Cal A}_{{\Cal V},
x_0}(\rho)}{V_{x_0}(\rho)}.
$$
Stoll showed that ${\Cal V}$ is algebraic if and only if $
\nu_{\infty} ({\Cal V})<\infty$. The following result generalizes
his result to analytic sets in curved manifold (we only succeeded
in codimension one case).

\medskip

{\it
 Let $M$ be a complete K\"ahler manifold
with nonnegative bisectional curvature. Let ${\Cal V}$ be a
analytic divisor of $M$. Then
 ${\Cal V}$ is defined by a `polynomial function' (holomorphic function of
  polynomial growth) if and only if
 $\nu_{\infty}({\Cal V})<\infty$.
 }

\medskip

In Section 4 we also obtained other results,  one of which
generalizes  the classical transcendental B\'ezout estimate for
codimension one analytic sets. (It has been known that the result
fails for the high codimension case by the famous example of
Cornalba and Shiffman \cite{CS}, even for the Euclidean case.) We
found the connection between the parabolic approach, especially
the Li-Yau-Hamilton estimate, and the classical Nevanlinna theory
 interesting and believe in that the parabolic approach should be the
 most nature/effective approach in extending sharp results from
 Euclidean spaces (linear) to the curved complex manifolds (nonlinear, in certain sense).

Finally, we discussed the relation between the Li-Yau-Hamilton
inequality proved in this paper and the previous computations in
\cite{FIN} on the reduced distance and the reduced volume modelled
on  Ricci expanders. In particular, we prove  a sharp lower bound
on the heat kernel for the time dependent case, which, in a sense,
is dual to Perelman's  monotonicity of the reduced volume. We also
explain how one can view this result, and more importantly,
Perelman's monotonicity of the reduced volume as a nonlinear
version of earlier work of Cheeger-Yau \cite{CY} and Li-Yau
\cite{LY} on the heat kernel estimates for heat
equations/Schr\"odinger equations. We also proved several local
monotonicity formulae for Ricci flow on the {\it forward reduced
volume} defined in \cite{FIN} and Perelman's  entropy, {\it
without any curvature sign assumption}. This again follows the
observation of Ecker in \cite{Ec2}, where he derived a localized
version of Huisken's monotonicity formula for the mean curvature
flow. As an application, we prove that the blow-down limit of a
so-called type III solution to K\"ahler-Ricci flow with bounded
nonnegative bisectional curvature must be a gradient expanding
soliton. This is, in certain sense, dual to Perelman's result on
ancient solutions.

Here is how we organize the paper. In Section 1 we prove the
interpolating version of the new Li-Yau-Hamilton estimate, which
in particular, includes (0.1). In Section 2 we derive some
monotonicity formulae and the heat kernel comparison theorem. In
Section 3 we show the interpolation between the new
Li-Yau-Hamilton inequality and  Perelman's entropy formula in
\cite{P}. In Section 4 we derive the localized version and prove
the manifold version of Stoll's theorem. In Section 5,  we discuss
the relation of the new inequality with the work of \cite{FIN},
formulate a new matrix Li-Yau-Hamilton expression on curvature
tensors for Ricci flow, and show a sharp heat kernel lower bound
estimate.

\medskip

{\it Acknowledgement}. The special thanks go to Professor B. Chow
as one can see that his contribution is crucial to a couple of
results in this paper. He  generously encouraged the author to
publish the results alone even though it should  really be  a
joint paper. The author would also like to thank Professors H.-D.
Cao, Tom Ilmanen, Peter Li and Jiaping Wang for helpful
discussions, Professors Klaus Ecker, Luen-Fai Tam, H. Wu for their
interests.

\input amstex
\documentstyle{amsppt}
\magnification=1200 \hsize=13.8cm \catcode`\@=11
\def\NoLogo{\let\logo@\empty}
\catcode`\@=\active \NoLogo

\def\heat{\lf(\frac{\p}{\p t}-\Delta\ri)}

\def \b {\beta}

\def\Ric{\text{Ric}}
\def\lf{\left}
\def\ri{\right}
\def\bbar{\bar \beta}
\def\a{\alpha}

\def\g{\gamma}
\def\e{\epsilon}
\def\p{\partial}
\def\delbar{{\bar\delta}}

\def\dbar{\bar\partial}

\def\C{\Bbb C}
\def\R{\Bbb R}
\def\tN{\tilde N}
\def\tY{\tilde Y}

\def\cN{\Cal N}
\def\ctN{{\tilde {\Cal N}}}

\def\cQ{\Cal Q}

\def\vp{\varphi}

\def\tN{\tilde N}

\def\dbar{\bar\partial}

\def\bb{{\bar\beta}}
\def\abb{{\alpha\bar\beta}}
\def\gbd{{\gamma\bar\delta}}

\def \D {\Delta}
\def\aint{\frac{\ \ }{\ \ }{\hskip -0.4cm}\int}
\documentstyle{amsppt}
\vsize=19.0 cm

\subheading{\S1 A new matrix Li-Yau-Hamilton inequality for
K\"ahler-Ricci flow}

\vskip .2cm

Let $M^m$ be a complete K\"ahler manifold of complex dimension
$m$. Let $(M, g(t))$ be a solution to K\"ahler-Ricci flow:
$$
\frac{\p}{\p t}g_{\abb}(x,t)=-R_{\abb}(x,t). \tag 1.1
$$
Here $R_{\abb}(x,t)$ is the Ricci tensor of the metric
$g_{\abb}(x,t)$. Let $u$ be a positive solution to the {\it
forward conjugate heat equation}:
$$
\heat u(x,t)={\Cal R}(x,t)u(x,t). \tag 1.2
$$

Here ${\Cal R}(x,t)$ is the scalar curvature. We shall prove the
following new matrix Li-Yau-Hamilton/differential Harnack
inequality.

\proclaim{Theorem 1.1} Let $(M, g(t))$ be a solution to {\rm
(1.1)} with nonnegative bisectional curvature. In the case that
$M$ is complete noncompact, assume further that the bisectional
curvature is bounded. Let $u$ be a positive solution to {\rm
(1.2)}. Then
$$
u_{\abb}+\frac{u}{t}g_{\abb}+ uR_{\abb} +u_{\a}V_{\bb} +u_{\bb}
V_{\a}+uV_{\a}V_{\bb}\ge 0 \tag 1.3
$$
for any $(1,0)$ vector field $V$.
\endproclaim

The assumption on $(M, g(t))$ having bounded nonnegative
bisectional curvature can be replaced by $(M, g(0))$ having
nonnegative bisectional curvature and $(M, g(t))$ with bounded
curvature, thanks to the result of Bando \cite{B}, Mok \cite{M2}
and Shi \cite{Sh2}. The same applies to other results of the
paper. We state the result under the stronger assumption just for
the simplicity.

Recall that in \cite{CN} the authors proved that for a fixed
K\"ahler metric $(M, g_{\abb}(x))$ with the nonnegative
bisectional curvature and the positive solution $u$ to the heat
equation $\heat u=0$ one has the matrix Li-Yau-Hamilton
inequality:
$$
u_{\abb}+\frac{u}{t}g_{\abb} +u_{\a}V_{\bb} +u_{\bb}
V_{\a}+uV_{\a}V_{\bb}\ge 0 \tag 1.4
$$
Therefore, one can think (1.3) is the nonlinear version of (1.4).
In fact we shall show that there exists a linear interpolation
between these two inequalities. The similar interpolation was
established by Ben Chow \cite{Ch1} originally for the Li-Yau's
gradient estimates for the heat equation and Hamilton's
differential Harnack for the Ricci flow in the case of surfaces.
We shall show such interpolation between the matrix differential
inequalities (1.3) and (1.4) for any dimensions. In \cite{N4} we
have shown another family of Li-Yau-Hamilton inequalities which
also serves as a higher dimensional generalization of Chow's
result.

Let us first set up the notations. For any $\epsilon >0$, we
consider the K\"ahler-Ricci flow:
$$
\frac{\p}{\p t} g_{\abb}(x,t)=-\e R_{\abb}(x,t). \tag 1.5
$$
Consider the positive solution $u$ to the parabolic equation:
$$
\heat u(x,t)=\e {\Cal R}(x,t)u(x,t). \tag 1.6
$$
We shall call (1.6) {\it forward conjugate heat equation}. We
shall prove the following general interpolation theorem.

\proclaim {Theorem 1.2} (Chow-Ni) Assume that $(M, g(t))$ has
nonnegative bisectional curvature satisfying {\rm (1.5)}. In the
case that  $M$ is complete noncompact, additionally assume that
the bisectional curvature is bounded. Let $u$ be a positive
solution to {\rm (1.6)}. Then
$$
u_{\abb}+\frac{u}{t}g_{\abb}+ \e uR_{\abb} +u_{\a}V_{\bb} +u_{\bb}
V_{\a}+uV_{\a}V_{\bb}\ge 0 \tag 1.7
$$
for any $(1,0)$ vector field $V$.
\endproclaim

\proclaim{Remark 1.1} Since Ben Chow helped with the proof at a
critical step it is well-justified to attribute the above theorem
as a joint result.
\endproclaim

It is easy to see that Theorem 1.2, serving an interpolation
between (1.3) and (1.4),  implies Theorem 1.1 and the earlier
result (1.4) for the linear heat equation. Therefore, we only need
to prove  Theorem 1.2. The proof  consists of the following
several lemmas.

\proclaim{Lemma 1.1}
$$
({\Cal R})_{\abb} =\D R_{\abb} +R_{\abb \gbd}R_{\delbar\g}-R_{\a
\bar{p}}R_{p\bar{\b}}. \tag 1.8
$$
Here $R_{\abb\gbd}$ is the bisectional curvature.
\endproclaim
\demo{Proof} The second Bianchi identity yields
$$
\split {\Cal R}_{\abb}&= (R_{\g\bar{\g}})_{,\abb}=R_{\a\bar{\g},\g\bar{\b}} \\
&= R_{\a\bar{\g}, \bar{\b}\g}-R_{p\bar{\g}\g\bar{\b}}R_{\a
\bar{p}}+R_{\a \bar{p}\g\bar{\b}}R_{p\bar{\g}}\\
&=R_{\abb, \bar{\g}\g}-R_{p\bar{\b}}R_{\a \bar{p}}+R_{\a
\bar{p}\g\bar{\b}}R_{p\bar{\g}}.
\endsplit
\tag 1.9
$$
The commutator formula gives that
$$
\split R_{\abb, \bar{\g}\g}&=R_{\abb,
\g\bar{\g}}-R_{p\bar{\b}\g\bar{\g}}R_{\a\bar{p}}+R_{\a\bar{p}\g\bar{\g}}R_{p\bar{\b}}\\
& =R_{\abb,
\g\bar{\g}}-R_{p\bar{\b}}R_{\a\bar{p}}+R_{\a\bar{p}}R_{p\bar{\b}}.
\endsplit \tag 1.10
$$
The lemma now follows from the definition $\D
R_{\abb}=\frac{1}{2}\lf(R_{\abb, \g\bar{\g}}+R_{\abb,
\bar{\g}\g}\right).$
\enddemo

\proclaim{Lemma 1.2} Let $u(x,t)$ be a solution to {\rm (1.6)}.
Then
$$
\heat
u_{\abb}=R_{\abb\gbd}u_{\bar{\g}\delta}-\frac{1}{2}\lf(R_{\a\bar{p}}u_{p\bar{\b}}+
R_{p\bar{\b}}u_{\a\bar{p}}\ri)+\e({\Cal R} u)_{\abb}.\tag 1.11
$$
\endproclaim
\demo{Proof} Differentiate (1.6) we have
$$
(u_t)_{\g\delbar}=R_{\beta \bar{\a}\g\delbar}u_{\a\bbar}+
g^{\a\bbar}u_{\a\bbar\g\delbar}+\e({\Cal R}u)_{\gbd}. \tag 1.12
$$
By definition $\Delta u_{\a\bbar}=
\frac{1}{2}\left(u_{\a\bbar,\g\bar{\g}}+u_{\a\bbar,
\bar{\g}\g}\right)$, in  normal coordinates at a point. We  need
to calculate the difference between the partial derivative
$u_{\a\bbar\g\delbar}$ and the covariant derivative
$u_{\a\bbar,\g\delbar}$. Direct computations show that, for normal
coordinates at a point
$$
u_{\g\delbar,\a\bbar}=u_{\g\delbar
\a\bbar}+u_{s\delbar}R_{\a\bbar\g\bar{s}} .\tag 1.13
$$
Using the fact that
$$
u_{\g\delbar, \a\bar{\a}}=u_{\g\delbar,
\bar{\a}\a}+R_{\g\bar{p}}u_{p\delbar} -R_{p\delbar}u_{\g\bar{p}}
\tag 1.14
$$
and (1.13) we have
$$
\split \Delta u_{\g\delbar}& =\frac{1}{2}(u_{\g\delbar,
\a\bar{\a}}+u_{\g\delbar,
\bar{\a}\a})\\
& = u_{\g\delbar
\a\bar{\a}}+\frac{1}{2}\left(R_{\g\bar{p}}u_{p\delbar}+
R_{p\delbar}u_{\g\bar{p}}\right).
\endsplit \tag 1.15
$$
Combining the above with (1.12), we conclude that $u_\abb$
satisfies (1.11).
\enddemo

The direct calculations give the following lemma.

\proclaim{Lemma 1.3} Let $(M, g(t))$ be a solution to {\rm (1.5)}.
Let $u(x,t)$ be a positive solution to {\rm (1.6)}. Then
$$
\split
 \heat \lf(\frac{u_{\a}u_{\bar{\b}}}{u}\ri)&= \e\frac{({\Cal
 R}u)_{\a}u_{\bar{\b}}}{u}+\e\frac{({\Cal
 R}u)_{\bar{\b}}u_{\a}}{u}-\e({\Cal R}u)\frac{u_{\a}
 u_{\bar{\b}}}{u^2}-\frac{u_{\a\bar{s}}u_{\bar{\b}s}+u_{\a
 s}u_{\bar{\b}\bar{s}}}{u}\\
 &\quad  -2\frac{u_{\a}
 u_{\bar{\b}}|u_s|^2}{u}-\frac{1}{2}\frac{R_{\a\bar{s}}u_{s}u_{\bar{\b}}
 +R_{s\bar{\b}}u_{\bar{s}}u_{\a}}{u}\\
 &\quad
 +\frac{u_{\a s}u_{\bar{s}}u_{\bar{\b}}+u_{\bar{\b}s}u_{\a}u_{\bar{s}}
 +u_{\a\bar{s}}u_{\bar{\b}}u_{s}+u_{\bar{\b}\bar{s}}u_{\a}u_{s}}{u^2}.
 \endsplit \tag 1.16
$$
$$
\heat \lf(\frac{u}{t}g_{\abb}\ri)(x,t)=\e\frac{{\Cal
R}u}{t}g_{\abb}-\frac{u}{t^2}g_{\abb}-\e\frac{u}{t}R_{\abb}. \tag
1.17
$$
$$
\split \heat \lf(\e u R_{\abb}\ri)& =\e^2 {\Cal R}uR_{\abb}+\e u
\lf(R_{\abb\gbd}R_{\bar{\g}\delta}-R_{\a\bar{p}}R_{p\bar{\b}}\ri)\\
&\quad  \e\lf(\nabla_s u\nabla_{\bar{s}}R_{\abb}+\nabla_{\bar{s}}
u \nabla_s
R_{\abb}\ri)\\
&\quad +(\e^2-\e)u\lf(
R_{\abb\gbd}R_{\bar{\g}\delta}-R_{\a\bar{p}}R_{p\bar{\b}}+\D
R_{\abb}\ri).
\endsplit \tag 1.18
$$
\endproclaim
\demo{Proof} In the derivation of (1.16), a commutator formula has
been used. In the derivation of (1.18) we have used
$$
\frac{\p}{\p t}R_{\abb}=\e\lf(
R_{\abb\gbd}R_{\bar{\g}\delta}-R_{\a\bar{p}}R_{p\bar{\b}}+\D
R_{\abb} \ri).
$$
\enddemo

We also need the following result which is an easy consequence of
Cao's  differential Harnack for the K\"ahler-Ricci flow \cite{C1}.

\proclaim{Lemma 1.4} Let $(M, g(t))$ be a complete solution to
{\rm (1.5)} with nonnegative bisectional curvature. In the case
that  $M$ is complete noncompact we further assume that the
bisectional curvature is bounded. Let $u(x,t)$ be a positive
solution to {\rm (1.6)}. Then
$$
\D R_{\abb}+R_{\abb\gbd}R_{\bar{\g}\delta}-\lf(\frac{\nabla_s
u}{\e u}\nabla_{\bar{s}}R_{\abb}+\frac{\nabla_{\bar{s}} u}{\e u}
\nabla_s R_{\abb}\ri)+R_{\abb\gbd}\frac{\nabla_{\bar{\g}}u}{\e
u}\frac{\nabla_{\delta}u}{\e u}+\frac{R_{\abb}}{\e t}\ge 0. \tag
1.19
$$
\endproclaim

Once we have Lemma 1.1--1.4  we can give the proof of Theorem 1.2.

\demo{Proof of Theorem 1.2} Let
$$
N_{\abb}=u_{\abb}+\frac{u}{t}g_{\abb}-\frac{u_{\a}
u_{\bar{\b}}}{u}
$$
and
$$
\tN_{\abb}=N_{\abb}+\e u R_{\abb}.
$$
By taking the minimizing vector in (1.7), one can see that Theorem
1.2 is equivalent to $\tN_{\abb}\ge 0$. Since the maximum
principle  needs some growth conditions for the noncompact
manifolds, we first prove the theorem for the case that $M$ is
compact.

\demo{Compact case} By Lemma 1.1--1.3, we have that
$$
\split \heat \tN_{\abb}
&=R_{\abb\gbd}u_{\bar{\g}\delta}-\frac{1}{2}\lf(R_{\a\bar{p}}u_{p\bar{\b}}+
R_{p\bar{\b}}u_{\a\bar{p}}\ri)+\e({\Cal R} u)_{\abb} \\
&\quad   -\e\frac{({\Cal
 R}u)_{\a}u_{\bar{\b}}}{u}-\e\frac{({\Cal
 R}u)_{\bar{\b}}u_{\a}}{u}+\e({\Cal R}u)\frac{u_{\a}
 u_{\bar{\b}}}{u^2}+\frac{u_{\a\bar{s}}u_{\bar{\b}s}+u_{\a
 s}u_{\bar{\b}\bar{s}}}{u}\\
 &\quad  +2\frac{u_{\a}
 u_{\bar{\b}}|u_s|^2}{u}+\frac{1}{2}\frac{R_{\a\bar{s}}u_{s}u_{\bar{\b}}
 +R_{s\bar{\b}}u_{\bar{s}}u_{\a}}{u}\\
 &\quad
 -\frac{u_{\a s}u_{\bar{s}}u_{\bar{\b}}+u_{\bar{\b}s}u_{\a}u_{\bar{s}}
 +u_{\a\bar{s}}u_{\bar{\b}}u_{s}+u_{\bar{\b}\bar{s}}u_{\a}u_{s}}{u^2}\\
 &\quad  +\e\frac{{\Cal
R}u}{t}g_{\abb}-\frac{u}{t^2}g_{\abb}-\e\frac{u}{t}R_{\abb}\\
&\quad  +\e^2 {\Cal R}uR_{\abb}+\e u
\lf(R_{\abb\gbd}R_{\bar{\g}\delta}-R_{\a\bar{p}}R_{p\bar{\b}}\ri)\\
&\quad  -\e\lf(\nabla_s u\nabla_{\bar{s}}R_{\abb}+\nabla_{\bar{s}}
u \nabla_s
R_{\abb}\ri)\\
&\quad +(\e^2-\e)u\lf(
R_{\abb\gbd}R_{\bar{\g}\delta}-R_{\a\bar{p}}R_{p\bar{\b}}+\D
R_{\abb}\ri)\\
&= R_{\abb\gbd}\tN_{\bar{\g}\delta}-\frac{1}{2}\lf(R_{\a
\bar{p}}\tN_{p\bar{\b}}+R_{p\bar{\b}}\tN_{\a
\bar{p}}\right)+\frac{1}{u}N_{\a\bar{p}}N_{p\bar{\b}}-\frac{2}{t}\tN_{\abb}\\
& \quad +\frac{1}{u}\lf(u_{\a p}-\frac{u_\a
u_p}{u}\ri)\lf(u_{\bar{p}\bar{\b}}-\frac{u_{\bar{p}}u_{\bar{\b}}}{u}\ri)+\e
{\Cal R} \tN_{\abb}-u\e^2 R_{\a \bar{p}}R_{p\bar{\b}}\\
&\quad +\e^2 u \tY_{\abb}
\endsplit \tag 1.20
$$
where
$$
\tY_{\abb}=\D
R_{\abb}+R_{\abb\gbd}R_{\bar{\g}\delta}-\lf(\frac{\nabla_s u}{\e
u}\nabla_{\bar{s}}R_{\abb}+\frac{\nabla_{\bar{s}} u}{\e u}
\nabla_s R_{\abb}\ri)+R_{\abb\gbd}\frac{\nabla_{\bar{\g}}u}{\e
u}\frac{\nabla_{\delta}u}{\e u}+\frac{R_{\abb}}{\e t}.
$$
By Lemma 1.4 and (1.20) we have that
$$
\split
 \heat \tN_{\abb} &= R_{\abb\gbd}\tN_{\bar{\g}\delta}-\frac{1}{2}\lf(R_{\a
\bar{p}}\tN_{p\bar{\b}}+R_{p\bar{\b}}\tN_{\a
\bar{p}}\right)+\lf(\e {\Cal R}-\frac{2}{t}\ri) \tN_{\abb}\\
&\quad +\frac{1}{2u}\tN_{\a \bar{p}}\lf(N_{p\bar{\b}}-\e u
R_{p\bar{\b}}\ri) +\frac{1}{2u}\lf(N_{\a\bar{p}}-\e u R_{\a
\bar{p}}\ri)\tN_{p \bar{\b}}\\
&\quad +\e^2u\tilde Y_{\abb}+\frac{1}{u}\lf(u_{\a p}-\frac{u_\a
u_p}{u}\ri)\lf(u_{\bar{p}\bar{\b}}-\frac{u_{\bar{p}}u_{\bar{\b}}}{u}\ri)\\
&\ge R_{\abb\gbd}\tN_{\bar{\g}\delta}-\frac{1}{2}\lf(R_{\a
\bar{p}}\tN_{p\bar{\b}}+R_{p\bar{\b}}\tN_{\a
\bar{p}}\right)+\lf(\e {\Cal R}-\frac{2}{t}\ri) \tN_{\abb}\\
&\quad +\frac{1}{2u}\tN_{\a \bar{p}}\lf(N_{p\bar{\b}}-\e u
R_{p\bar{\b}}\ri) +\frac{1}{2u}\lf(N_{\a\bar{p}}-\e u R_{\a
\bar{p}}\ri)\tN_{p \bar{\b}}.
\endsplit \tag 1.21
$$
Using the observation that the right hand side of (1.21) satisfies
the null-vector condition of the tensor maximum principle of
Hamilton \cite{H2}, we have proved the case $M$ being compact.

\enddemo

\demo{Noncompact case} First recall  the fundamental derivative
estimate of Shi.   For $g_{\abb}(x,t)$,  a solution to (1.5) on
$M\times [0, T]$ with bounded (in space-time)  nonnegative
bisectional curvature,  there exist $A_k>0$ such that
$$
\| \nabla^k R_{\abb\gbd}\|^2\le \frac{A_k}{t^k} \tag 1.22
$$
on $M\times [0, T]$. For our consideration we only need to prove
it for $T$ small. The estimate (1.22) is proved  in \cite{Sh1}.
For the sake of simplicity,  we will show the matrix differential
Harnack inequality Theorem 1.2 for the case $\e=1$  under the
above assumption (1.22). In fact, what needed is (1.22) for $k\le
2$. We also need the perturbation trick from \cite{NT1} (see also
\cite{CN}). Namely we first shift $t$ by $2\delta$, where $\delta$
is a small positive number. After the shifting we can have
estimates  on $u$,  $|\nabla u|$ and $|u_{\abb}|$. The goal is to
show that there exits a $b>0$ such that
$$
\int_\delta ^T \int_M
\exp(-b(r_0^2(x)+1))\left(\frac{1}{u}+\frac{|\nabla
u|^2}{u}+|u_{\abb}|^2\right)\, d\mu \, dt <\infty. \tag 1.23
$$
Here $r_0(x)$ is the distance to $x$ from a fixed point $o\in M$
with respect to the initial metric.  We need the estimate (1.23)
to apply the tensor maximum principle from  \cite{N5} (see also
\cite{NT2} for the original time-independent version).

In order to get control on $u$ (or $\frac{1}{u}$) first we need
the following Harnack inequality of Guenther \cite{Gu}.

\proclaim{Theorem 1.3} (Guenther) Let $(M, g(t))$ be a solution to
Ricci flow satisfying {\rm (1.22)}. Let $u$ be a positive solution
to the {\it forward conugate heat equation} $\heat u={\Cal R}u$.
Then for sufficient small $T$ (only depending  on $A_k$) there
exist $\a, B_1 >0$ only depending on $A_k$  such that
$$
u(x_1, t_1)\le u(x_2, t_2)
\left(\frac{t_2}{t_1}\right)^{\frac{n}{2}}\exp(\frac{r^2(x_1, x_2,
t_1)}{\a (t_2-t_1)}+B_1(t_2-t_1)). \tag 1.24
$$
for any $T\ge t_2>t_1>0$. Here $r(x_1, x_2, t_1)$ denotes the the
distance between $x_1$ and $x_2$ with respect to the metric at
$t=t_1$.
\endproclaim
The result above was proved in \cite{Gu} through a gradient
estimate of Li-Yau type on compact manifold. Since one can apply
the localization techniques as in \cite{LY, page 161} (see also
\cite{NT1, page 647}, \cite{Sh1} ) one can easily generalize the
gradient estimate, thus the above Harnack estimate, proved in
\cite{Gu} to complete noncompact case. From (1.24), one can deduce
that for small $\delta$ there exists a constant $b_2, B_2>0$,
where $b_2 =b_2(A_k, \delta)$ and $B_2=B_2(u(o, \frac{\delta}{4}),
u(o, T-\frac{\delta}{4}), A_k, \delta)$ such that
$$
\left(\frac{1}{u}+u\right)(x,t)\le B_2\exp (b_2 (r^2_0(x)+1)) \tag
1.25
$$
for $(x,t)\in M\times[\frac{\delta}{2}, T-\frac{\delta}{2}]$.
\enddemo
Observe that we have the following two equations.
$$
\heat u^2={\Cal R} u^2-|\nabla u|^2 \tag 1.26
$$
and
$$
\heat |\nabla
u|^2=-|u_{\abb}|^2-|u_{\a\beta}|^2+\langle\nabla({\Cal R}u),
\nabla u\rangle+\langle\nabla  u, \nabla({\Cal R}u)\rangle . \tag
1.27
$$
The desired estimate (1.23) follows from (1.26), (1.27) by the
argument through integration by parts in \cite{CN, Lemma 3.1}.
Once we have established the estimate (1.23) one can apply the
perturbation argument as in \cite{NT1} (see also \cite{CN})
together with the tensor maximum principle in \cite{NT2, Theorem
2.1} or \cite{N5, Theorem 2.1} to conclude the proof of the matrix
Li-Yau-Hamilton estimate (1.7) for the complete noncompact case.
Note that Theorem 1.3 only provide the estimate for short time.
But we can iterate the argument to prove the result for all time.
An alternative is that once one has the upper bound estimate at
some earlier time one can also make use of the heat kernel
estimate in \cite{N5} for the time dependent heat operator (and
uniqueness of the positive solution) to get estimates of the
positive solution for the later time. We also should remark that
the matrix Li-Yau-Hamilton (1.7) gives a sharp Harnack (which is
more precise, compared with Theorem 1.3). Please see Corollary 2.1
in the next section. However, we do need the rough estimate in the
proof to apply the tensor maximum principle. \enddemo

\proclaim{Corollary 1.1} Let $u(x,t)$ be a positive solution to
{\rm (1.6)}. Then
$$
\D \log u +\e {\Cal R}+\frac{m}{t}\ge 0. \tag 1.28
$$
If the equality holds for some $(x_0,t_0)$ with $t_0>0$, then
 $(M, g(t))$ is an expanding K\"ahler-Ricci
soliton for the case $\e
>0$,  and $(M, g)$ is isometric to $\C^m$ for the case $\e =0$.
\endproclaim
\demo{Proof} Let
$$
Q=\D \log u +\e {\Cal R}+\frac{m}{t}.
$$
and $\cQ =g^{\abb}\tY_{\abb}$. Since $\tilde N_{\abb}\ge 0$ and
$Q=g^{\abb}\tilde N_{\abb}$ we have that $\tilde N_{\abb}(x_0,
t_0)= 0$. By the strong maximum principle we know that
$\tN_{\abb}\equiv 0$ for all $t< t_0$. $\tilde N_{\abb}\equiv 0$
is nothing but the equation in the definition of the gradient
expanding soliton.

For the case $\e=0$, we apply the same line of argument. In this
case we have $Q\equiv 0$ since $\tilde N_{\abb}=N_{\abb}\equiv 0$,
for $t\le t_0$. Now from the equation
$$\split
0&=\heat t^2Q =t^2 R_{\abb}(\log u)_{\bar{\a}}(\log
u)_{\beta}+\frac{t^2}{u}\lf(u_{\a p}-\frac{u_\a
u_p}{u}\ri)\lf(u_{\bar{p}\bar{\a}}-\frac{u_{\bar{p}}u_{\bar{\a}}}{u}\ri)\\
&\quad +\frac{t^2}{u}N_{\a \bar{p}}N_{p\bar{\a}}\\
&\ge 0
\endsplit
$$
we have that $(\log u)_{\a\b}\equiv 0$ too. From the definition of
$N_{\abb}$ we then have
$$(\log u)_{\abb}=\frac{1}{t} g_{\abb}
$$
and
$$
(\log u)_{\a\beta}=0
$$
from which it is easy to see that $(M, g)$ is flat and in fact
isometric to $\C^m$ , since that the curvature (see, for example,
\cite{KM, page 117}) can be written as
$$
R_{\abb\gbd}=-\frac{\p^4 (tf)}{\p z^{\a}\p z^{\bar{\b}}\p z^{\g}\p
z^{\bar{\delta}}}+g^{p\bar{q}}\lf(\frac{\p^3 (tf)}{\p z^{\bar{q}}
\p z^{\a} \p z^{\g} }\ri)\lf(\frac{\p^3 (tf)}{\p z^{p} \p
z^{\bar{\b}}\p z^{\bar{\delta}}}\ri)
$$ and that $f=\log u$ (defined to be) is a convex function.
\enddemo

\proclaim{Remark 1.2} The case $\e=0$ case in the above Corollary
1.1  is just the original Li-Yau's estimate \cite{LY}, which holds
with nonnegativity of the Ricci curvature. The equality case
implies the manifold is $\R^n$ is implicit in the proofs of
\cite{LY} and proved explicitly in \cite{N3} for Riemannian
manifolds with nonnegative Ricci curvature. It would be
interesting to see if Corollary 1.1 is true, for $\e>0$ case, for
complete solutions to the K\"ahler-Ricci/Ricci flow without
assumptions on the sign of the curvature.
\endproclaim

When the manifold is compact, since we known that, by passing to
its universal cover, it is products of $\C^k$ with compact
Hermitian symmetric spaces. Without the loss of generality we can
assume that the first Chern class $c_1(M)$ is  a positive multiple
of the K\"ahler class. Then the K\"ahler-Ricci flow will have
singularity, say at $t=1$, and the normalized flow has long time
existence. The rescaling is given by $\hat{g}=\frac{1}{1-t}g$ and
the reparametrization is given by $s=-\log (1-t)$. Therefore,
Theorem 1.2 has the following equivalent form.

\proclaim{Theorem 1.2'} Let $M$ be a compact K\"ahler manifold as
above. Let $g(x,t)$ be a solution to {\rm (1.1)},  and let
$u(x,t)$ be a positive solution to {\rm (1.2)}. Assume that
$\hat{g}(x,s)$ is the solution to the normalized K\"ahler Ricci
flow. Then
$$
(\log u)_{\abb}+\hat{R}_{\abb}+\frac{1}{e^s-1}\hat{g}_{\abb}\ge
0.\tag 1.29
$$
Here $\hat{R}_{\abb}$ is the Ricci tensor of $\hat{g}$.
\endproclaim

\input amstex
\documentstyle{amsppt}
\magnification=1200 \hsize=13.8cm \catcode`\@=11
\def\NoLogo{\let\logo@\empty}
\catcode`\@=\active \NoLogo

\def\heat{\lf(\frac{\p}{\p t}-\Delta\ri)}

\def \b {\beta}

\def\Ric{\text{Ric}}
\def\lf{\left}
\def\ri{\right}
\def\bbar{\bar \beta}
\def\a{\alpha}

\def\g{\gamma}
\def\e{\epsilon}
\def\p{\partial}
\def\delbar{{\bar\delta}}

\def\dbar{\bar\partial}

\def\C{\Bbb C}
\def\R{\Bbb R}
\def\P{\Bbb P}
\def\E{\Bbb E}

\def\vp{\varphi}

\def\tN{\tilde N}
\def\cN{\Cal N}
\def\ctN{\tilde{\Cal N}}

\def\dbar{\bar\partial}

\def\bb{{\bar\beta}}
\def\abb{{\alpha\bar\beta}}
\def\gbd{{\gamma\bar\delta}}

\def \D {\Delta}
\def\aint{\frac{\ \ }{\ \ }{\hskip -0.4cm}\int}
\documentstyle{amsppt}
\vsize=19.0 cm

\subheading{\S2 Monotonicity formulae}

\vskip .2cm

In this section we derive some monotonicity formulae out of the
matrix Li-Yau-Hamilton  estimate, as well as its trace,  proved in
last section. In order to make the argument unified for both
cases,  with and without Ricci flow, we work with the
interpolation version (1.7). Let $(M, g(t))$ be the solution to
the Ricci flow (1.5) and let $u(x,t)$ be the positive solution to
(1.6).

The first result is the monotonicity of the partition function
(also called Nash's entropy in \cite{FIN}) defined by
$$\ctN(g, u, t)=-\int_M u\log u\, dv-m\log (\pi t)-m. \tag 2.1$$
Simple computation shows that
$$
\frac{d \ctN}{d t}=\int_M \left(-\D \log u-\epsilon{\Cal
R}-\frac{m}{t}\right)u\, d\mu_t\le 0. \tag 2.2
$$
Here $d\mu_t$ is the volume element of $g(t)$. The following
result is a direct consequence of Corollary 1.1.

\proclaim{Proposition 2.1} Let $(M, g(t))$ be a solution to {\rm
(1.5)} with nonnegative bisectional curvature.  Let $H(x,y, t)$ be
the fundamental solution to {\rm (1.6)}. Then
$$\frac{d}{d t}\ctN(g, u, t)\le 0$$
for any positive solution $u$ to {\rm (1.6)} and
$$
-\int_M H\log H\, d\mu_t -m\log(\pi  t) \le m
$$
with the equality holds for some positive $t$ if and only if the
manifold is an expanding gradient soliton (isometric to $\C^m$ in
the case $\e =0$).
\endproclaim

For the fixed  Riemannian metric case the above was proved earlier
in \cite{N2} for the manifold with only nonnegative Ricci
curvature. We believe that the result should hold for Ricci flow
even without assumptions. But at this moment we can only prove it
for K\"ahler-Ricci flow through the matrix Li-Yau-Hamilton
inequality, which assumes the nonnegativity of the bisectional
curvature. The result above gives a characterization of expanding
solitons using the partition function. The similar formulation for
the shrinking solitons also works for the partition functions
related to Perelman's entropy formula \cite{P}.  Namely, consider
the backward Ricci flow $\frac{\p}{\p \tau}g_{ij}=2R_{ij}$ on
$M\times [0, \tau_0]$, where  $M$ is a Riemannian manifold of real
dimension $n$. Let $u(x, \tau)$ be a solution to the {\it backward
adjoint heat equation} $$\lf(\frac{\p}{\p \tau} -\D +{\Cal
R}\ri)u(x,\tau)=0. \tag 1.6'$$ Similarly one can define
$$
\ctN(g, u, \tau)=-\int_M u\log u\, d\mu_\tau -\frac{n}{2}\log
(4\pi \tau)-\frac{n}{2}.\tag 2.1'
$$
In Proposition 1.2 of \cite{P},  Perelman proved that
$$
\frac{d \ctN}{d \tau}=\int_M \lf(-\D \log u+{\Cal
R}-\frac{n}{2\tau}\ri)\, d\mu_\tau \le 0. \tag 2.2'
$$
The dual version of Proposition 2.1 states as follows.

\proclaim{Proposition 2.1'} Let $H(x,y,\tau)$ be the fundamental
solution to the adjoint heat equation. Then
$$
-\int_M H\log H\, d\mu_\tau -\frac{n}{2}\log (4\pi \tau) \le
\frac{n}{2}
$$
with equality holds (or in {\rm (2.2')}) for some positive $\tau$
if and only if $(M, g_{ij}(\tau))$ is a gradient shrinking
soliton.
\endproclaim
\demo{Proof} Since we do not  have Corollary 1.1 in this case we
need other arguments. In fact, tracing the equality case of the
proof of Proposition 1.2 of \cite{P} we have that
$R_{ij}-\nabla_i\nabla_j \log H$ is diagonal. Namely
$R_{ij}-\nabla_i\nabla_j \log H=\frac{{\Cal R} -\D\log
H}{n}g_{ij}$. On the other hand, the equality in (2.2') further
implies that $R_{ij}-\nabla_i\nabla_j \log
H=\frac{1}{2\tau}g_{ij}$.
\enddemo
 In \cite{N2}, the relation between the value of $\ctN(t)$ (as well as the entropy functional
 ${\Cal W}$) as $t\to
 \infty$ and  the asymptotic volume ratio at infinity (also called the `cone angle')
 was proved for the linear heat equation. The similar
 relation,  between  the $\kappa$-constant in the
 $\kappa$-non-collapsing of volume defined in \cite{P} and  the
 large time limit of the partition function $\ctN(g,u,\tau)$ (as well he entropy functional
 ${\Cal W}$) should also be true for
 ancient solutions to  Ricci flow.

Another application of Theorem 1.2 is an entropy monotonicity for
the  ancient solutions. Let $(M, g(t))$ be an ancient solution to
(1.5) and  let $u$ be a positive solution to (1.6), both defined
on $M\times (-\infty, 0]$. Theorem 1.2 implies that
$$
(\log u)_{\abb}+R_{\abb}\ge 0.
$$
From this one can easily see that
$$\cN(g, u, t):=-\int_M u\log u\, d\mu_t
$$
is monotone non-increasing.

 If $M$ is compact, one can
obtain such $u$ by taking limit of solutions with initial data at
$t=-i$ as in \cite{FIN} (where the immortal solution is studied).
More precisely, let $u_i(x, t)$ be a solution to (1.6) with $u(x,
-i)=\frac{1}{V(-i)}$. Letting $i\to \infty$, one can extract a
limit $u_\infty>0$ (since $\int_M u_i\equiv 1$ and $u_i>0$) which
is defined on $(-\infty, 0]$. In this case $\int_M u_\infty\, dv
=1$. Applying the Jensen's inequality one has that
$$
\cN(g, u_\infty, t)\le \log V(t).
$$
The right hand side is another monotone non-increasing quantity
along the flow. When $M$ is complete  noncompact one can obtain
such a $u$ similarly by solving $u_i(x,t)$ with initial condition
at $t=-i$ and  anchoring $u_i(o, -1)=1$, where $o\in M$ is a fixed
point. In fact one can even get $u_i$ integrable by taking it to
be scalar multiple of the fundamental solution (with initial data
being the  delta function at $t=-i$).

The second application is to derive the heat kernel comparison
theorem and Huisken type \cite{Hu, E1-2} monotonicity formula for
the analytic subvarieties in $M$. We start with the case $\e=0$
since the results seem to be  more useful at this moment.

\proclaim{Theorem  2.1} Let $M$ be a complete K\"ahler manifold
with nonnegative bisectional curvature. Let $H(x,y, t)$ be the
fundamental solution of the heat equation. Let ${\Cal V}\subset M$
be any complex subvariety of dimension $s$. Let $K_{\Cal V}(x,y,
t)$ be the fundamental solution of heat equation on ${\Cal V}$.
Then \roster \item"{(i)}"
$$
K_{\Cal V}(x,y, t)\le (\pi t)^{m-s}H(x,y, t), \text{ for any }x,\,
y \in {\Cal V}. \tag 2.3
$$
If the equality holds, then  ${\Cal V}$ is totally geodesic.
Furthermore if $\tilde M$ is the universal cover of $M$ with
covering map $\pi$ and $\tilde {\Cal V}=\pi^{-1}({\Cal V})$, then
$\tilde M=\tilde M_1\times \C^k$ for some K\"ahler manifold
$\tilde M_1$ which does not contain any Euclidean factors, with
$k\ge m-s$. Moreover $\tilde {\Cal V}=\tilde M_1\times \C^l$ with
$l< k$. \item"{(ii)}"
$$
\frac{d}{d t}\int_{\Cal V} (\pi t)^{m-s}H(x,y, t) \, dA_{\Cal
V}(y) \ge 0, \text { for any }\, x\in M. \tag 2.4
$$
Similarly, if the equality holds for some $x\in M$ at some
positive time $t$, then   $\tilde M=\tilde M_1\times \C^k$ with
$k\ge m-s$.
\endroster
\endproclaim

\demo{Proof}  For any smooth point $y\in {\Cal V}$, choose a
complex coordinate $(z_1\, \cdots, z_m)$ such that $(z_1, \cdots,
z_s)$ is the coordinates for ${\Cal V}$. Let $i,\, j, \, k,
\cdots$ denote the coordinate on ${\Cal V}$ and $a, b, c, \cdots$
denote the coordinates in the normal directions. Then we compute
$$
\split \left(\Delta^{(y)}_{\Cal V}-\frac{\p}{\p t}\right)\lf( (\pi
t)^{m-s}H(x,y, t)\ri) & =\lf(\D_{M}
-g^{a\bar{b}}\nabla_{a}\nabla_{\bar{b}} -\frac{\p}{\p
t}\ri)\lf(t^{m-s} H(x,y, t)\ri)
\\
& =  (\pi t)^{m-s}\lf(\Delta_M
H-H_t-g^{a\bar{b}}\nabla_a\nabla_{\bar{b}}H -\frac{m-s}{t}H\ri)\\
&= (\pi t)^{m-s}\lf(-\nabla_a\nabla_{\bar{b}}H+\frac{\nabla_a H
\nabla_{\bar{b}}H}{H}-\frac{1}{t}Hg_{a\bar{b}}\ri)g^{a\bar{b}}\\
&\quad   -(\pi t)^{m-s}\frac{|\nabla ^{\perp} H|^2}{H}.
\endsplit \tag 2.5
$$
By Theorem 1.2 we have that
$$
\left(\Delta^{(y)}_{\Cal V}-\frac{\p}{\p t}\right)\lf(( \pi
t)^{m-s}H(x,y, t)\ri)\le 0.
$$
Noticing that for $x\in {\Cal V}$, $\lim_{t\to 0}( \pi
t)^{m-s}H(x,y, t)|_{\Cal V}=\delta_x(y)$. This proves that (2.3)
by the maximum principle,  (2.4) by the integrating (2.5) on
${\Cal V}$. (For the case ${\Cal V}$ is singular, one can refer to
\cite{LT} for the justification on the validity of the integration
by parts.)

In the case when the  equality holds in (2.3) we have that $(\pi
t)^{m-s}H(x,y, t)$ satisfies the heat equation.  Hence equality
holds in (2.5) for any $x, y\in {\Cal V}$. This implies that
$$
|\nabla^{\perp} \log H|^2 \equiv 0, \text{  for }x,\, y \in {\Cal
V}, \, t>0,  \tag 2.6
$$
which implies that $<\nabla r^2(x,y), \nu>=0$ for any $x, y\in
{\Cal V}$, and any normal direction $\nu$. Here we have used the
fact that $\lim_{t\to 0}-t\log H(x,y, t)=r^2(x,y)$. See
\cite{CLY1}, for example. This implies that for any smooth point
$x\in {\Cal V}$, any minimizing geodesic starting from $x$ lies
totally inside ${\Cal V}$. More precisely, let $\gamma_{v}(s)$ be
a short minimizing geodesic emitting from $x$ with
$\gamma'_{v}(0)=v$ such that $v\in T_{x}{\Cal V}$. Denote by $h$
be the distance function from ${\Cal V}$. Then
$$
\split \frac{d }{d s}h(\gamma_v(s))&=2Re \langle\nabla h, \gamma'_v(s)\rangle\\
&=\frac{1}{s}2Re\langle \nabla h, s\gamma'_v(s)\rangle\\
&=\frac{1}{s}Re\langle \nabla h, \nabla r^2(x,
\gamma_{v}(s))\rangle =0.
\endsplit
$$
Here $\nabla h=\nabla^\a h\frac{\p}{\p z^\a}$ and
$\gamma'_{v}(s)=\frac{d z^{\a}(s)}{ d s}\frac{\p}{\p z^\a}$.
Therefore, ${\Cal V}$ is totally geodesic. This in particular
shows that there is no singular points in ${\Cal V}$, which rules
out the possibility that ${\Cal V}$ is union of several totally
geodesic submanifolds with singular intersections (in which case
the heat kernel comparison (2.3) has strictly inequality).
Moreover, by lifting the computation to the universal cover we can
assume that $M$ is simply-connected. Then the manifold $M$ splits
by Theorem 0.1 of \cite{NT2}, more precisely, Theorem 2.1 and
Corollary 2.1 of \cite{NT2}, since $ N_{\abb}=\nabla
\a\nabla_{\bar{\b}}H-\frac{\nabla_\a H
\nabla_{\bar{\b}}H}{H}+\frac{1}{t}Hg_{\abb}\ge 0$ and its null
space is at least of $m-s$ dimension, by (2.5), for any $x,\,
y\in{\Cal V}$ and $t>0$. Notice that $N_{\abb}$  does not
satisfies the linear Lichnerowicz heat equation. However the
inequality (1.21), satisfied by $N_{\abb}$ (since $\e=0$),  is
enough for the argument in the proof of Corollary 2.1 of
\cite{NT2}. Then $M=M_1\times M_2$ such that the tangent space of
$M_2$ consists of the null space of $N_{\abb}$. The  factor $M_2$
is isometric to $\C^{m-s}$ follows from the same argument in the
proof of Corollary 1.1, since
$-\nabla_a\nabla_{\bar{b}}H+\frac{\nabla_a H
\nabla_{\bar{b}}H}{H}-\frac{1}{t}Hg_{a\bar{b}}\equiv 0$.
 (One can also use Corollary 1.3 of \cite{N2, page 331}.) More precisely, if we
 write a point $x\in M$ as
$x=(x_1, x_2)$ according to the splitting, we can write the heat
kernel $H(x,y,t)=H_1(x_1, y_1, t) H_2(x_2, y_2, t)$. Then on $M_2$
we have that $(\log H_2(x_2, y_2))_{\abb}+\frac{1}{t}g_{\abb}=0$
by the definition of the splitting. Therefore one can apply
Corollary 1.1 to conclude that $M_2=\C^k$.
 If (2.4)  holds equality for some $x\in
 M$, it implies that the right hand side of (2.5) is zero. Then
 the argument as above also applies. Note that we may not have ${\Cal
 V}$ totally geodesic since $x$ may be a singular point.
\enddemo

\proclaim{Remark 2.1} 1) In the case ${\Cal V}$ is not smooth, the
existence on $K(x,y,t)$  was justified in the work of Li and Tian
\cite{LT}. They also obtained a similar upper bound estimate as
{\rm (2.3)} for the special case $M=\P^m$  with the Fubini-Study
metric using very different method. Their method produces better
upper bound for the special case $M=\P^m$.   Our estimate here
works for general K\"ahler manifolds with nonnegative bisectional
curvature.

2) The similar heat kernel comparison was first proved in
\cite{CLY2} for minimal submanifolds in space forms by quite
different method. The result in Theorem 2.1 is more general than
\cite{CLY2} in the sense that it holds for any K\"ahler manifolds
with nonnegative sectional curvature in stead of space forms.
However  it is also more restrictive since  it only applies to
analytic subvarieties.

3) In the part (ii) of Theorem 2.1, the manifold $\tilde M_1$ may
not be $\tilde {\Cal V}$. This could happen, for example, in the
case $M=\C^m$ and ${\Cal V}$ is a union of two hyper-planes and
$x$ lies on the intersection subvariety. If we further assume that
equality holds for all smooth points $x\in {\Cal V}$, we do have
the same conclusion as part (i).

4) The monotonicity {\rm (2.4)} is enough for  applications in
\cite{N4}. Namely, one can prove the comparison results on the
dimensions of polynomial growth holomorphic function spaces
obtained in \cite{N4},  using {\rm (2.4)} in stead of the other
Li-Yau-Hamilton inequality proved therein. In fact one can derive
the results in \cite{N4} through the following corollary, which is
a special case of  Theorem 2.1 (or Corollary 2.2) of \cite{N4}. It
is however enough for the applications considered in \cite{N4}.
\endproclaim

\proclaim{Corollary 2.1} Let $f$ be a holomorphic function. Let
${\Cal V}=Z(f)$, the zero locus of $f$. Denote
$$
w(x,t)=\int_M H(x,y,t) \D \log |f|^2(y)\, d\mu(y). \tag 2.7
$$
Then
$$
tw(x,t)=(\pi t)\int_{\Cal V}  H(x,y, t)\, dA_{\Cal V}(y). \tag 2.8
$$
Moreover
$$
\frac{\p}{\p t} \left(t w(x,t)\right)\ge 0. \tag 2.9
$$
If the equality holds for some point $x\in M$ and some positive
time $t$, then  the universal cover (of $M$) $\tilde M$ splits at
least a factor of $\C$.
\endproclaim
\demo{Proof} Notice that $\D =g^{\abb}\frac{\p}{\p z^{\a}\p
z^{\bar{\b}}}$, which differs a factor of $4$ from \cite{N4}. Let
$f$ be a holomorphic function defined on $M$. Here we do require
that $f$ does not growth too fast. For example, it will be
sufficient if $f$ is of polynomial growth or is of finite order in
the sense of Hadamard (see (3.1) of \cite{N4} for precise
definition). We can write
$$v(x,t)=\int_M H(x,y,t) \log |f|^2(y)\, d\mu(y)$$ as a solution
to the heat equation with initial value $\log|f|^2(y)$.  The
requirement on $f$ is to make such representation formula of
$v(x,t)$ meaningful. Then $w(x,t)$ is defined to be $\frac{\p}{\p
t} v(x,t)$ as in \cite{N4}. By the definition of  $w(x,t)$ in
(2.7) it is easy to see that $w(x,t)$ is also a solution to the
heat equation with initial data given by the measure $\D \log
|f|^2$. It is easy to see that $w(x,t)=\frac{\p}{\p t} v(x,t)$. By
the Poincar\'e-Lelong formula we know that
$$
(\pi t)\int_M H(x,y,t)\lf(\frac{\sqrt{-1}}{2\pi}\p\dbar \log
|f|^2(y)\ri)\wedge \frac{\omega^{m-1}}{(m-1)!}=(\pi t)\int_{\Cal
V} H(x,y,t)\, dA_{\Cal V}(y).
$$
On the other hand, the direct calculation shows that the left hand
side of the above equation is equal to
$$
\split
 t\int_M H(x,y,t)\D \log |f|^2(y) \,
\frac{\omega^m}{m!}&=t\int_M H(x,y,t) \D \log |f|^2(y)\, d\mu(y)\\
&=tw(x,t).
\endsplit
$$
This proves the first statement of the corollary. The monotonicity
(2.9)  is just a special case  (ii) of Theorem 2.1. The proof
above also shows that $w(x,t)$ has the same meaning as in Lemma
3.1 and Theorem 3.1, Theorem 4.1 of \cite{N4}.
\enddemo

The case with Ricci flow, namely $\e=1$,  can be formulated
similarly. In order to do so we have to explain some notations.
For a time-dependent metrics deformed by the K\"ahler-Ricci flow
equation (1.1),
 we call $H(x,y, t, t_0)$ a  fundamental solution to the heat
equation, if for any $u(x,t_0)=\vp(x)$, where $\vp(x)$ is a
compact-supported smooth function, the solution $u(x,t)$ to the
heat equation $\heat u(x,t)=0$ with initial condition
$u(x,t_0)=\vp(x)$ is given by the formula
$$
u(x,t)=\int_M H(x,y,t, t_0)\vp(y)\, d\mu(y).
$$
It is easy to check that $H(x,y,t, t_0)$ must satisfies the {\it
forward conjugate heat equation} (1.2). When the meaning is clear
in the context (mostly $t_0=0$) we just simply write $H(x,y,t,
t_0)$ as $H(x,y, t)$. Adapting the notation used in the proof of
Theorem 2.1, the restriction of the K\"ahler-Ricci flow (1.5) on
${\Cal V}$ reads as $\frac{\p }{\p
t}g_{i\bar{j}}(x,t)=-R_{i\bar{j}}(x,t)$. We denoted
$g^{i\bar{j}}R_{i\bar{j}}$ by ${\Cal R}_{\Cal V}$. It is easy to
see that ${\Cal R}_{\Cal V}$ is well-defined. Therefore the
fundamental solution to the heat equation on ${\Cal V}$, $K_{\Cal
V}(x,y,t)$ satisfies
$$
\heat v(x,t)={\Cal R}_{\Cal V}(x,t) v(x,t). \tag 2.10
$$
We have the following comparison and monotonicity result.

\proclaim{Theorem 2.1'} Let $M$ be a complete K\"ahler manifold
with bounded nonnegative bisectional curvature. Let $H(x,y,t)$ be
a fundamental solution to the heat equation on $M$. Let ${\Cal V}$
be a complex subvariety of $M$ of dimension $s$. Let $K_{\Cal
V}(x,y,t)$ be the fundamental solution to the heat equation (with
respect to the induced metrics) on ${\Cal V}$. Then we have  {\rm
(2.3)} and {\rm (2.4)}. Moreover, the equality (for positive $t$),
in either cases, implies that the universal cover (of $M$) $\tilde
M$
 has the splitting $\tilde M=\tilde M_1\times \E^k$, where $\E^k$ is an gradient
 expanding K\"ahler-Ricci
soliton of dimension $k\ge m-s$.
\endproclaim

\proclaim{Remark 2.2} One can think {\rm (2.4)} as a dual version
of Perelman's monotonicity of the reduced volume since the reduced
volume of \cite{P, Section 7} is, in a sense,  a`weighted volume'
of $M$ (with weight being the heat kernel of a `potentially
infinity dimensional manifold' restricted to $M$, as explained in
Section 6 of \cite{P}), while here the monotonicity is on the
`weighted volume' of complex submanifolds with weight being the
heat kernel of $M$ restricted to the submanifold. The reduced
volume monotonicity of Perelman  has important applications in the
study of Ricci flow. We expect that {\rm (2.4)} will have some
applications in understanding the K\"ahler-Ricci flow and its
effect on the complex geometry of analytic subvarieties.
\endproclaim

Taking the trace of the matrix estimate in Theorem 1.2 and
integrating  along the space-time path as in \cite{LY}, we can
have the following Harnack estimates for the positive solutions to
the {\it forward conjugate heat equation}. This gives a sharp
version of the previous rough estimate of Guenther in \cite{Gu}.

\proclaim{Corollary 2.1} Let $(M, g(t))$ be a solution to Ricci
flow {\rm (1.1)} and $u(x,t)$ be a positive solution to {\rm
(1.2)}. Then
$$
\frac{|\nabla u|^2}{u^2}-\frac{u_t}{u}+\frac{m}{t}\ge 0
$$
and  for any $t_2>t_1$,
$$
u(x_2, t_2)t_2^m \ge u(x_1, t_1)t_1^m
\exp\lf(-\inf_{\g}\int_{t_1}^{t_2}|\g'(t)|^2\, dt\ri) \tag 2.11
$$
Here $\g(t)$ is a path with $\g(t_1)=x_1$ and $\g(t_2)=x_2$.
\endproclaim

Note that we do not have the factor $4$ due to our choice of $\D$
and the gradient $|\nabla f|^2=g^{\abb}f_\a f_{\bar{\b}}$. Also
$$
|\gamma'(t)|^2=g_{\abb}\frac{dz^\a}{d t}\frac{d z^{\bar{\b}}}{d
t}.
$$
We list this consequence here since it implies the monotonicity of
$t^m u(x,t)$.

\input amstex
\documentstyle{amsppt}
\magnification=1200 \hsize=13.8cm \catcode`\@=11
\def\NoLogo{\let\logo@\empty}
\catcode`\@=\active \NoLogo

\def\heat{\lf(\frac{\p}{\p t}-\Delta\ri)}

\def \b {\beta}

\def\Ric{\text{Ric}}
\def\lf{\left}
\def\ri{\right}
\def\bbar{\bar \beta}
\def\a{\alpha}

\def\g{\gamma}
\def\e{\epsilon}
\def\p{\partial}
\def\delbar{{\bar\delta}}

\def\dbar{\bar\partial}

\def\C{\Bbb C}
\def\R{\Bbb R}
\def\tZ{\tilde Z}

\def\cN{\Cal N}
\def\ctN{{\tilde {\Cal N}}}

\def\cQ{\Cal Q}

\def\vp{\varphi}

\def\tN{\tilde N}

\def\dbar{\bar\partial}

\def\bb{{\bar\beta}}
\def\abb{{\alpha\bar\beta}}
\def\gbd{{\gamma\bar\delta}}

\def \D {\Delta}
\def\aint{\frac{\ \ }{\ \ }{\hskip -0.4cm}\int}
\documentstyle{amsppt}
\vsize=19.0 cm

\subheading{\S3 Interpolation between Perelman's entropy formula
and the new Li-Yau-Hamilton inequality}

\vskip .2cm

The purpose of this section is two-folded. First we  give a
different proof of Theorem 1.2.  The second purpose is to show
that the by-product of this different proof  also implies
Perelman's monotonicity of entropy, as well as the energy. The
computation in this section has its real version. See \cite{Ch2}
and the forth-coming book \cite{CLN}. The main computation is
summarized in equation (3.7) below, which we call a  {\it
pre-Li-Yau-Hamilton} equality. The equation (3.7) can also be
viewed as a matrix version of Perelman's entropy monotonicity
formula. In a sense, one can think that  the matrix
Li-Yau-Hamilton inequality proved in Section 1 is dual to
Perelman's  entropy formula in \cite{P, Section 3}.

Consider the K\"ahler-Ricci flow:
$$
\frac{\p}{\p \tau} g_{\abb} =\epsilon R_{\abb} \tag 3.1
$$
where $\epsilon$ is a parameter and the {\it conjugate heat
equation}:
$$
\lf( \frac{\p}{\p \tau} -\D +\epsilon {\Cal R}\ri) u(x,\tau) =0.
\tag 3.2
$$
When $\e<0$, (3.1) is a forward Ricci flow equation and (3.2)
become {\it forward conjugate heat equation}. The equations look
different from those in Section 1 since in this section the case
of $\e<0$ corresponds to the forward Ricci flow and the case of
$\e>0$ corresponds to the backward Ricci flow. For example $\e=1$
is exactly the setting for Perelman's entropy and energy
monotonicity. Notice that (3.2) becomes the {\it backward
conjugate heat equation} for $\e=1$. For the positive solution
$u(x,\tau)$ we define the $(1,1)$ tensor $Z_{\abb}$ by
$$
Z_{\abb}=-(\log u)_{\abb} +\epsilon R_{\abb}.
$$
Let $\D_L$ denote the Lichnerowicz Laplacian on $(1,1)$ tensors,
which is defined by
$$
\D_L\eta_{\abb} =\D\eta_{\abb} +R_{\abb
\gbd}\eta_{\bar{\g}\delta}-\frac{1}{2}\lf( R_{\a \bar{p}}\eta_{p
\bar{\b}}+\eta_{\a \bar{p}}R_{p \bar{\b}}\ri)
$$
for any Hermitian symmetric (1.1) tensor $\eta_{\abb}$. It is
known, see for example \cite{C1}, that
$$
\lf(\frac{\p}{\p \tau}-\D_L\ri)R_{\abb}=-(1+\epsilon) \D_L
R_{\abb}.
$$
The direct calculation as in \cite{NT1, Lemma 2.1} shows that

\proclaim{Lemma 3.1} For any $C^2$-function $f(x,\tau)$
$$
\lf(\frac{\p}{\p \tau} -\D_L\ri) f_{\abb} =\lf[\lf(\frac{\p}{\p
\tau}-\D \ri)f\ri]_{\abb}.\tag 3.3
$$
\endproclaim

\proclaim{Remark 3.1} The similar result as the above lemma   hold
for Ricci flow on Riemannian manifolds. Please see \cite{CLN} for
details.
\endproclaim

Now with the help of Lemma 3.1 and 1.1 we can calculate $
\lf(\frac{\p}{\p \tau}-\D_L\ri)Z_{\abb} $ as follows using the
equation $\lf(\frac{\p}{\p \tau}-\D_L\ri)(\log u)=-\epsilon{\Cal
R}+|\nabla \log  u|^2$.
$$
\split \lf(\frac{\p}{\p \tau}-\D_L\ri)Z_{\abb}&=\lf(\epsilon{\Cal
R}-|\nabla \log u|^2\ri)_{\abb}-\epsilon(1+\epsilon)\D_L R_{\abb}\\
&= \epsilon({\Cal R})_{\abb}-\lf(g^{\gbd}(\log u)_{\g}(\log
u)_{\bar{\delta}}\ri)_{\abb}-\epsilon(1+\epsilon)\D_L R_{\abb}\\
&=-R_{\abb\gbd}(\log u)_{\bar{\g}}(\log u)_{\delta}-(\log u)_{\a
\g}(\log u)_{\bar{\g}\bar{\b}}-(\log
u)_{\a \bar{\g}}(\log u)_{\g \bar{\b}}\\
&\quad \quad -\lf[(\log u)_{\abb}\ri]_{\g} (\log u)_{\bar{\g}}
-\lf[(\log u)_{\abb}\ri]_{\bar\g} (\log u)_{\g}-\epsilon^2\D_L
R_{\abb}.
\endsplit
\tag 3.4
$$
Hence
$$
\split \lf(\frac{\p}{\p \tau}-\D_L\ri)Z_{\abb}&=-\epsilon^2\lf(\D
R_{\abb}+R_{\abb\gbd}
R_{\bar{\g}\delta}+\nabla_{\g}R_{\abb}(\frac{1}{\epsilon}\nabla_{\bar\g}\log
u)+\nabla_{\bar\g}R_{\abb}(\frac{1}{\epsilon}\nabla_{\g}\log u) \ri.\\
&\quad  +\lf.R_{\abb\gbd}(\frac{1}{\epsilon}\nabla_{\bar{\g}}\log
u)
(\frac{1}{\epsilon} \nabla_{\delta}\log u)\ri)\\
&\quad  +\epsilon^2R_{\a \bar{\g}}R_{\g\bar{\b}}-(\log u)_{\a
\g}(\log u)_{\bar{\g}\bar{\b}}-(\log u)_{\a \bar{\g}}(\log u)_{\g
\bar{\b}}
\\
&\quad  +\nabla_{\g}(Z_{\abb})\nabla_{\bar{\g}}\log
u+\nabla_{\bar{\g}}(Z_{\abb})\nabla_{\g}\log u.
\endsplit\tag 3.5
$$
Regrouping terms yields
$$
\split \lf(\frac{\p}{\p \tau}-\D_L\ri)Z_{\abb}&= -\epsilon^2\lf(\D
R_{\abb}+R_{\abb\gbd}
R_{\bar{\g}\delta}+\nabla_{\g}R_{\abb}(\frac{1}{\epsilon}\nabla_{\bar\g}\log
u)+\nabla_{\bar\g}R_{\abb}(\frac{1}{\epsilon}\nabla_{\g}\log u) \ri.\\
&\quad  +\lf.R_{\abb\gbd}(\frac{1}{\epsilon}\nabla_{\bar{\g}}\log
u)
(\frac{1}{\epsilon} \nabla_{\delta}\log u)\ri)\\
&\quad  -(\log u)_{\a \g}(\log u)_{\bar{\g}\bar{\b}}
+\nabla_{\g}(Z_{\abb})\nabla_{\bar{\g}}\log
u+\nabla_{\bar{\g}}(Z_{\abb})\nabla_{\g}\log u\\
&\quad  +\frac{1}{2}Z_{\a\bar{\g}}\lf(\epsilon R_{\g
\bar{\b}}+(\log u)_{\g \bar{\b}}\ri)+\frac{1}{2}\lf(\epsilon
R_{\a\bar{\g}}+(\log u)_{\a\bar{\g}}\ri)Z_{\g \bar{\b}}.
\endsplit\tag 3.6
$$
Let
$$
\tZ_{\abb}=Z_{\abb}-\frac{1}{\tau}g_{\abb}
$$
and
$$
\split
 Y_{\abb}&=\D R_{\abb}+R_{\abb\gbd}
R_{\bar{\g}\delta}+\nabla_{\g}R_{\abb}(\frac{1}{\epsilon}\nabla_{\bar\g}\log
u)+\nabla_{\bar\g}R_{\abb}(\frac{1}{\epsilon}\nabla_{\g}\log u)\\
&\quad  +R_{\abb\gbd}(\frac{1}{\epsilon}\nabla_{\bar{\g}}\log u)
(\frac{1}{\epsilon} \nabla_{\delta}\log u).
\endsplit
$$
Notice  that $\tilde Y_{\abb}$, defined after (1.20)  in Section
1, is related to $Y_{\abb}$ above through the equation $\tilde
Y_{\abb}=Y_{\abb}-\frac{R_{\abb}}{\e\tau}$ (remember that $-\e$
here corresponding to $\e$ in Section 1). From (3.6), we can
derive the equation for $\tZ_{\abb}$ as follows. \proclaim{Lemma
3.2} (Chow-Ni)
$$
\split \lf(\frac{\p}{\p \tau}-\D_L\ri)\tZ_{\abb}&=
\lf(\frac{\p}{\p
\tau}-\D_L\ri)Z_{\abb}+\frac{1}{\tau^2}g_{\abb}-\frac{1}{\tau}\epsilon
R_{\abb}\\
&= -\lf(\e^2Y_{\abb}-\frac{\e}{\tau}R_{\abb}\ri)
-(\log u)_{\a \g}(\log u)_{\bar{\g}\bar{\b}}\\
&\quad  +\nabla_{\g}(\tZ_{\abb})\nabla_{\bar{\g}}\log
u+\nabla_{\bar{\g}}(\tZ_{\abb})\nabla_{\g}\log u\\
&\quad  +\frac{1}{2}\tZ_{\a\bar{\g}}\lf(\e R_{\g \bar{\b}}+(\log
u)_{\g \bar{\b}}\ri)+\frac{1}{2}\lf(\e R_{\a\bar{\g}}+(\log
u)_{\a\bar{\g}}\ri)\tZ_{\g \bar{\b}}-\frac{1}{\tau}\tZ_{\abb}.
\endsplit
\tag 3.7
$$
\endproclaim
Notice that $\tZ_{\abb}=\frac{1}{u}\tilde N_{\abb}$. With some
labor one can check that (1.21) and (3.7) are equivalent. Namely
one can derive one from the other, keeping in mind that $-\e$ here
corresponds $\e$ in Section 1. One can also write (3.7) as
$$
\split \lf(\frac{\p}{\p \tau}-\D_L\ri)\tZ_{\abb}&=
\lf(\frac{\p}{\p
\tau}-\D_L\ri)Z_{\abb}+\frac{1}{\tau^2}g_{\abb}-\frac{1}{\tau}\epsilon
R_{\abb}\\
&= -\lf(\e^2Y_{\abb}-\frac{\e}{\tau}R_{\abb}\ri)
-(\log u)_{\a \g}(\log u)_{\bar{\g}\bar{\b}}\\
&\quad  +\nabla_{\g}(\tZ_{\abb})\nabla_{\bar{\g}}\log
u+\nabla_{\bar{\g}}(\tZ_{\abb})\nabla_{\g}\log u\\
&\quad +\frac{1}{2}\tZ_{\a\bar{\g}}\lf(\e R_{\g \bar{\b}}+(\log
u)_{\g \bar{\b}}-\frac{1}{\tau}g_{\g
\bar{\beta}}\ri)\\
&\quad  +\frac{1}{2}\lf(\e R_{\a\bar{\g}}+(\log
u)_{\a\bar{\g}}-\frac{1}{\tau}g_{\a \bar{\g}}\ri)\tZ_{\g
\bar{\beta}}.
\endsplit
\tag 3.7'
$$

For $\e<0$,  applying Lemma 1.4, the fundamental result of H.-D.
Cao, we know that $\e^2Y_{\abb}-\frac{\e}{\tau}R_{\abb}\ge 0$
under the assumption that $M$ is a complete K\"ahler manifold with
bounded nonnegative holomorphic bisectional curvature. Hence the
tensor maximum principle and (3.7) implies that $\tZ_{\abb}\le 0$,
which is equivalent to the statement of Theorem 1.2. Namely, (3.7)
does lead to  another proof  of Theorem 2.1.

Even though the computation (3.7) and (1.21) are essentially
equivalent, (3.7) has the advantage that when $\e=1$ it also
implies Perelman's energy/entropy monotonicity formulae. The
following is a more detailed computation of this claim. Let
$f=-\log u$, $Z=g^{\abb}Z_{\abb}$. The tracing of (3.6) gives
$$
\split \lf(\frac{\p}{\p \tau}-\D\ri)Z &=-\epsilon
R_{\bar{\a}\beta}Z_{\abb}-\epsilon^2g^{\abb}Y_{\abb}-
(f)_{\a\g}(f)_{\bar{\g}\bar{\a}}-\nabla_{\g}Z \nabla_{\bar{\g}}f-\nabla_{\bar{\g}}Z\nabla_{\g}f \\
&\quad +Z_{\abb}(\e R_{\b\bar{\a}}-f_{\b\bar{\a}})
\endsplit \tag 3.8
$$
and
$$
g^{\abb}Y_{\abb}=\D {\Cal
R}+R_{\abb}R_{\bar{\a}\b}-\nabla_{\g}{\Cal
R}(\frac{1}{\e}\nabla_{\bar{\g}}f)-(\frac{1}{\e}\nabla_{\g}f)\nabla_{\bar{\g}}{\Cal
R}+R_{\abb}(\frac{1}{\e}f_{\bar{\a}})(\frac{1}{\e}f_{\b}).
$$
The following observation of Ben Chow is also useful.
\proclaim{Lemma 3.3} (Chow) In the case $\e=1$, we have that
$$
\int_M \left(g^{\abb}Y_{\abb}\right) u\, d\mu =\int_M
\lf(R_{\abb}(R_{\bar{\a}\b}+f_{\bar{\a}\b})\ri)u\, d\mu. \tag 3.9
$$
\endproclaim
\demo{Proof}  Follows from the integration by parts and the second
Bianchi identity ${\Cal R}_{\g}=R_{\g \bar{\a}, \a}$.
\enddemo

\proclaim{Remark 3.2} Please refer to \cite{Ch2} and \cite{CLN}
for the Riemannian version of the above identity.
\endproclaim
Recall the definition  of energy ${\Cal F}$.
$$
{\Cal F}(g,u,\tau)=\int_M \left(\frac{|\nabla u|^2}{u}+{\Cal R}
u\right)\, d\mu.
$$
Then (3.8) (with $\e=1$)  and  Lemma 3.3 implies the the following
result. \proclaim{Proposition 3.1} (Perelman)
$$
\frac{d}{d \tau}{\Cal F}(g,u,\tau)=-\int_M
\left(|R_{\abb}+f_{\abb}|^2+|f_{\a\b}|^2\ri)u\, d\mu. \tag 3.10
$$
\endproclaim
Note that this is nothing but the energy monotonicity formula of
Perelman in Section 1 of \cite{P}. We show that it follows from
(3.8). \demo{Proof of Proposition 3.1} Let $\e=1$ in (3.8) we have
that
$$
\lf(\frac{\p}{\p \tau}-\D\ri)Z =-g^{\abb}Y_{\abb}-
(f)_{\a\g}(f)_{\bar{\g}\bar{\a}}-\nabla_{\g}Z \nabla_{\bar{\g}}f
-\nabla_{\bar{\g}}Z\nabla_{\g}f
 -Z_{\abb}f_{\b\bar{\a}}. \tag 3.8'
$$
Now the result follows from direct computation on $\frac{d}{d
\tau}{\Cal F}(g,u,\tau)= \frac{d}{d \tau}\int_M Z u\, d\mu$, by
applying Lemma 3.3.

\enddemo

Similarly, if we trace (3.7) and denote $\tZ=g^{\abb}\tZ_{\abb}$,
we have that
$$
\split \lf(\frac{\p}{\p \tau}-\D\ri)\tZ &=-\epsilon
R_{\bar{\a}\beta}\tZ_{\abb}-\epsilon^2g^{\abb}Y_{\abb}+\frac{\e}{\tau}{\Cal
R}-
(f)_{\a\g}(f)_{\bar{\g}\bar{\a}}-\nabla_{\g}\tZ \nabla_{\bar{\g}}f-\nabla_{\bar{\g}}\tZ\nabla_{\g}f \\
&\quad +\tZ_{\abb}(\e
R_{\b\bar{\a}}-f_{\b\bar{\a}})-\frac{1}{\tau}\tZ.
\endsplit \tag 3.11
$$
For $\e=1$, integration by parts as before gives
$$
\frac{d}{d \tau}\int_M \tZ u\, d\mu_{\tau}=-\int_M
\left(|\tZ_{\abb}|^2+|f_{\a\b}|^2\ri)u\,
d\mu_\tau-\frac{2}{\tau}\int_M \tZ u \, d\mu_\tau. \tag 3.12
$$
The above equation is equivalent to  Perelman's entropy
monotonicity formula due to the following consideration. Let
$$
\ctN(g, u, \tau):=-\int_M u\log u\, d\mu_\tau -m\log (\pi \tau)-m.
$$
Then Perelman's  entropy
$$
{\Cal W}(g, u, \tau)=\int_M\left[\tau(2\D \bar f-|\nabla \bar
f|^2+{\Cal R}+\bar f-2m\right]u\, d\mu_\tau,
$$
where $\bar f=-\log u-m\log(\pi \tau)$, can be expressed as
$$
{\Cal W}(g, u, \tau)=\frac{d}{d\tau} (\tau \ctN)=\tau \int_M \tZ
u\, d\mu_\tau+\ctN . \tag 3.13
$$
Therefore
$$
\split \frac{d}{d\tau}{\Cal W}&=\tau \frac{d}{d\tau}\int_M \tZ u\,
d\mu_\tau +2\int_M \tZ u\, d\mu_\tau\\
&=-\tau \int_M \left(|\tZ_{\abb}|^2+|f_{\a\b}|^2\ri)u\, d\mu_\tau
\endsplit \tag 3.14
$$
which is nothing but the entropy formula of Perelman in \cite{P,
Section 3}.

As pointed out in \cite{P}, there exists a statistical mechanics
analogy of Perelman's entropy. If we identify the quantities above
with the notation of \cite{Ev2, Chapter I and VII}, $\tau$ is the
temperature; $-\ctN$ defined above  is the $\log$ of the
distribution function in \cite{Ev2}; $-{\Cal W}$ is the entropy
$S$ in \cite{Ev2}; $\tau\ctN $ is the free energy $F$; $\frac{\p
\ctN}{\p ( \frac{1}{\tau})}$ is the energy $E$ in \cite{Ev2}
(which is nonnegative by Proposition 1.2 of \cite{P}) and the
first equation of (3.13) is just the well-known equation
$S=-\frac{\p F}{\p \tau}$ from thermodynamics. The negation of the
right hand side of (3.14), measuring the deviation from an
shrinking soliton, is called the {\it heat capacity} in
\cite{Ev2}. The entropy formula (3.14) implies the concavity of
$S$ in $E$, one of the defining properties for the entropy. This
analogy also holds for the solution to the heat equation with
respect to a fixed metric Riemannian metric with nonnegative Ricci
curvature \cite{N3}.

\proclaim{Remark 3.3}
  For the case of $\e=0$,
the similar computation as above gives the monotonicity formula in
\cite{N3}. The strange thing is that one can not get nice
monotonicity formula in the case $\e \ne -1$ or $0$. Namely, the
interpolation formula {\rm (3.6)/(3.7)} does not give nice
interpolation for energy/entropy after integration on $M$. One can
view {\rm (3.6)/(3.7)} as a matrix version of the energy/entropy
monotonicity formula for $\e>0$. This partially answers one of the
questions raised in the end of \cite{N3} (still not satisfactory
though). On the other hand, when $\e<0$, one can not get entropy
monotonicity out of {\rm (3.7)}. Instead we have  a pointwise
Li-Yau-Hamilton inequality. However, for $\e =0$, both entropy and
the differential Harnack follows from {\rm (3.6)} (See \cite{N3}).
The above discussion indicates that our new matrix Li-Yau-Hamilton
inequality is dual to the entropy monotonicity of Perelman in some
sense and there may perhaps  be certain profound duality behind
the scene.

Another puzzling point is that so far we have not been able to
verify a matrix Li-Yau-Hamilton inequality analogue to Theorem 1.2
for the Ricci flow on Riemannian manifolds, even though the above
computation {\rm (3.7)} holds for $\e>0$ for the Ricci flow on
Riemannian manifolds (which  was done in \cite{Ch2}). In short,
the interpolation between positive and negative $\e$ by now only
works in K\"ahler category. The validity of the Riemannian case is
pending on the verification of a new matrix Li-Yau-Hamilton
estimate similar to Hamilton's famous work \cite{H1}. Please see
Remark 5.2 for further details.
\endproclaim

\input amstex
\documentstyle{amsppt}
\magnification=1200 \hsize=13.8cm \catcode`\@=11
\def\NoLogo{\let\logo@\empty}
\catcode`\@=\active \NoLogo

\def\heat{\lf(\frac{\p}{\p t}-\Delta\ri)}

\def \b {\beta}

\def\Ric{\text{Ric}}
\def\lf{\left}
\def\ri{\right}
\def\bbar{\bar \beta}
\def\a{\alpha}

\def\g{\gamma}
\def\e{\epsilon}
\def\p{\partial}
\def\delbar{{\bar\delta}}

\def\dbar{\bar\partial}

\def\C{\Bbb C}
\def\R{\Bbb R}
\def\P{\Bbb P}
\def\E{\Bbb E}

\def\vp{\varphi}
\def\test{\vp_{(x_0, t_0), \rho}}
\def\testt{\vp_{(x_0, t_1), \rho}}

\def\tN{\tilde N}
\def\cN{\Cal N}
\def\ctN{\tilde{\Cal N}}

\def\dbar{\bar\partial}

\def\bb{{\bar\beta}}
\def\abb{{\alpha\bar\beta}}
\def\gbd{{\gamma\bar\delta}}

\def \D {\Delta}
\def\aint{\frac{\ \ }{\ \ }{\hskip -0.4cm}\int}
\documentstyle{amsppt}
\vsize=19.0 cm

\subheading{\S4 A local monotonicity formula and its applications}

\vskip .2cm

 This section is inspired by the work of
Ecker in \cite{E1-2}. Let $M$ be a complete K\"ahler manifold with
nonnegative bisectional curvature (unless specified otherwise). In
this section we study the localization of the previous established
monotonicity (in Section 2) for a fixed K\"ahler metric. (We leave
the K\"ahler-Ricci flow case to a later discussion.)  In
\cite{E1}, the localized  monotonicity formula is proved for mean
curvature flow in Euclidean spaces. Since we are dealing with
curved spaces here,  we need  some extra ingredients, which
includes  Theorem 2.1 in Section 2, the complex Hessian comparison
theorem on distance functions proved in \cite{LW} and \cite{CN}
recently and the well-known heat kernel estimates of Li-Yau on
complete Riemannian manifolds with nonnegative Ricci curvature,
which states that
$$
\frac{C^{-1}(n)}{V_x(\sqrt{t})}\exp\left(-\frac{r_x^2(y)}{3t}\right)\le
H(x,y,t)\le
\frac{C(n)}{V_x(\sqrt{t})}\exp\left(-\frac{r_x^2(y)}{5t}\right)
$$
for some $C(n)>0$, where $n$ is the real dimension of the manifold
considered. As applications we prove an elliptic `monotonicity
principle' for complex subvarieties in $M$. It can then be applied
to prove a manifold version of Stoll's theorem.

Before we prove a  localized version of Theorem 2.1, we need to
introduce some functions (notations). Let ${\Cal V}$ be an
analytic subvariety (of $M$) of complex dimension $s$. For the
simplicity of the notation, Let $H(x,y, \tau)$ be the fundamental
solution to the heat equation $\left(\frac{\p}{\p
\tau}-\D\right)u(x,\tau)=0$. When $x, y\in {\Cal V}$, we denote
$(\pi \tau)^{m-s}H(x,y, \tau)$ by $H_{\Cal V}(x,y, \tau)$. For any
fixed $(x_0, t_0)$ with $x_0\in {\Cal V}$,  we denote $H_{\Cal
V}(x_0, y, t_0-t)$ by $ \hat H_{(x_0, t_0), {\Cal V}}(y, t)$. For
any $\rho>0$,  we also introduce a cut-off function
$$
\vp_{(x_0, t_0), \rho}(y,
t)=\left(1-\frac{r^2_{x_0}(y)+s(t-t_0)}{\rho^2}\right)_{+} \tag
4.1
$$
where $f_{+}(x)=\max(f, 0)$ for any function $f$, $r_{x_0}(y)$ is
the distance function (of $M$) from $x_0$ to $y$. It is easy to
see that $\test$ is supported in $
B_{x_0}(\sqrt{\rho^2-s(t-t_0)})$.

The following simple lemma is useful.

\proclaim{Lemma 4.1} On ${\Cal V}$,
$$\left(\frac{\p}{\p
t}-\D_{\Cal V}\right)\test (y, t) \le 0. \tag 4.2$$ Here $\D_{\Cal
V}$ denotes the Laplacian operator with respect to the induced
K\"ahler metric on ${\Cal V}$ (strictly speaking only regular part
of ${\Cal V}$).
\endproclaim
\demo{Proof} For any $y\in {\Cal V}$, choose a complex coordinate
$(z_1\, \cdots, z_m)$ as in the proof of Theorem 2.1. Namely
$z_\a=0$ on ${\Cal V}$ for any $\a>s$. We also use the index
convention as in the proof of Theorem 2.1. Namely, $1\le i, j, k,
\cdots \le s$ and $s+1\le a,b,c,\cdots \le m$. Direct computation
shows that
$$
\split \left(\frac{\p}{\p t}-\D_{\Cal
V}\right)\left(s(t-t_0)+r^2_{x_0}(y)\right)&=
s-g^{i\bar{j}}\left(r^2_{i\bar{j}}\right)\\
&\ge s-g^{i\bar{j}}g_{i\bar{j}}\\
&\ge 0.
\endsplit
$$
Here we have used the Hessian comparison theorem on the distance
functions proved in \cite{LW} (see also  \cite{CN, Corollary
1.1}).
\enddemo

The following is a localized version of part (ii) of Theorem 2.1.
\proclaim{Proposition 4.1} Let
$$
E_{{\Cal V}, x_0, t_0, t_1}(t)=\int_{\Cal V} \testt(y, t)\,
 \hat H_{(x_0, t_0), {\Cal V}}(y, t)\, dA_{\Cal V}. \tag 4.3
$$When in the right context we also it briefly denote by $E_{\Cal V}$.
Then
$$
\frac{d}{d t}E_{\Cal V}(t)\le -\int_{\Cal V}|\nabla^{\perp}\log
 \hat H_{(x_0, t_0), {\Cal V}}|^2 \testt  \hat H_{(x_0, t_0), {\Cal V}}\,
dA_{\Cal V}. \tag 4.4
$$
\endproclaim
\demo{Proof} The computation (2.5) implies that
$$
\left(\frac{\p}{\p t}+\D_{\Cal V}\right) \hat H_{(x_0, t_0), {\Cal
V}}\le -|\nabla^{\perp}\log  \hat H_{(x_0, t_0), {\Cal V}}|^2 \hat
H_{(x_0, t_0), {\Cal V}}. \tag 4.5
$$
By Lemma 4.1 we have that
$$
\split
 \frac{d}{d t}E_{\Cal V}(t)&=\int_{\Cal V} \left(\frac{\p}{\p t}
 \testt\right)  \hat H_{(x_0, t_0), {\Cal V}}+ \testt \left(\frac{\p}{\p t}
   \hat H_{(x_0, t_0), {\Cal V}}\right)\, dA_{\Cal V}\\
&=\int_{\Cal V}\left(\left(\frac{\p}{\p t}-\D_{\Cal
V}\right)\testt \right) \hat H_{(x_0, t_0), {\Cal V}}+\testt
\left(\left(\frac{\p}{\p
t}+\D_{\Cal V}\right) \hat H_{(x_0, t_0), {\Cal V}}\right)\\
&\quad + \int_{\Cal V}(\D_{\Cal V}\testt) \hat H_{(x_0, t_0),
{\Cal V}}-\testt
\left(\D_{\Cal V}  \hat H_{(x_0, t_0), {\Cal V}}\right)\, dA_{\Cal V}\\
&\le -\int_{\Cal V}|\nabla^{\perp}\log  \hat H_{(x_0, t_0), {\Cal
V}}|^2 \testt  \hat H_{(x_0, t_0), {\Cal V}}\, dA_{\Cal V}.
\endsplit
$$
Here we have used the observation that
$$
\int_{\p({\Cal V}\cap B_{x_0}(\sqrt{\rho^2+s(t_0-t)}))} \hat
H_{(x_0, t_0), {\Cal V}}\langle \nabla_{\Cal V}\testt, \nu\rangle
\, dS \le 0
$$
where $\nu$ is the unit out-normal of $\p({\Cal V}\cap
B_{x_0}(\sqrt{\rho^2+s(t_1-t)}))$ in ${\Cal V}$ and $dS$ is the
area integral of $\p({\Cal V}\cap
B_{x_0}(\sqrt{\rho^2+s(t_1-t)}))$.
\enddemo
We denote the $2s$-dimensional Hausdorff measure of set ${\Cal
V}\cap B_{x_0}(\rho)$ by ${\Cal A}_{{\Cal V}, x_0}(\rho)$. As a
consequence we have the following elliptic `monotonicity
principle' (can be viewed a Bishop-Lelong lemma on manifolds).

\proclaim{Corollary 4.1} (Monotonicity principle) Let
$\delta(s)=\frac{1}{\sqrt{2+4s}}$. There exists $C=C(m, s)$ such
that for any $\rho'\in (0, \delta(s)\rho)$
$$
\frac{{\Cal A}_{{\Cal V}, x_0}(\rho')
(\rho')^{2(m-s)}}{V_{x_0}(\rho')}\le C(m,s) \frac{{\Cal A}_{{\Cal
V}, x_0}(\rho) \rho^{2(m-s)}}{V_{x_0}(\rho)}.\tag 4.6
$$
\endproclaim
\demo{Proof} The proof follows essentially the argument of
Proposition 3.5 in \cite{E2}. Applying Proposition 4.1 with
$t_1=t_0-\delta^2(s)\rho^2$ and $t_0$ replaced by $t_0+{\rho'}^2$.
Using the heat kernel upper bound of Li-Yau we have that
$$
\split
 H_{(x_0, t_0+{\rho'}^2), {\Cal V}}(y,
t_0-\delta^2(s)\rho^2)&\le \left(\pi (\delta^2(s)
\rho^2+{\rho'}^2)\right)^{m-s}
\frac{C(m)}{V_{x_0}(\sqrt{{\rho'}^2+\delta^2(s)\rho^2})}\\
&\le C(m,s)\frac{ \left(\pi \rho^2\right)^{m-s}}{V_{x_0}(\delta
(s)\rho)}\\
&\le C(m,s)\frac{ \left(\pi \rho^2\right)^{m-s}}{V_{x_0}(\rho)}.
\endsplit \tag 4.7
$$
Here we have used the Bishop volume comparison in the last
inequality. Notice that $\testt(x, t_0-\delta^2(s)\rho^2)\le 1$
and supported inside $B_{x_o}(\rho)$. Hence (4.7) implies that
$$
E_{\Cal V}(t_0-\delta^2(s)\rho^2)\le C(m,s) \frac{{\Cal A}_{{\Cal
V}, x_0}(\rho) \rho^{2(m-s)}}{V_{x_0}(\rho)}. \tag 4.8
$$
By Proposition 4.1 we know that
$$
E_{\Cal V}(t_0-\delta^2(s)\rho^2)\ge E_{\Cal V}(t_0-{\rho'}^2).
\tag 4.9
$$
On the other hand, by Li-Yau's heat kernel lower bound we also
have that, for all $y\in B_{x_0}({\rho'})$,
$$
\split
 H_{(x_0, t_0+{\rho'}^2), {\Cal V}}(y,
t_0-{\rho'}^2)&\ge \left(\pi (2 {\rho'}^2)\right)^{m-s}
\frac{C(m)}{V_{x_0}(\sqrt{2}\rho')}\exp\left(-\frac{r^2_{x_0}(y)}{6{\rho'}^2}\right)\\
&\ge  C(m,s)\frac{ \left(\pi
{\rho'}^2\right)^{m-s}}{V_{x_0}(\rho')}.
\endsplit \tag 4.10
$$
Again we have used the Bishop volume comparison theorem. Notice
that for $y\in B_{x_0}({\rho}')$
$$
\split
 \testt (y, t_0-{\rho'}^2) &\ge 1-
\frac{{\rho'}^2+s(\delta^2(s)\rho^2+{\rho'}^2)}{\rho^2}\\
&\ge 1-\delta^2(s)(1+2s)\\
&= \frac{1}{2}.
\endsplit \tag 4.11
$$
Combining (4.10), (4.11) we have that
$$
E_{\Cal V}(t_0-{\rho'}^2)\ge C(m,s)\frac{{\rho'}^{2(m-s)} {\Cal
A}_{{\Cal V}, x_0}({\rho'})}{V_{x_0}(\rho')}. \tag 4.12
$$
Combining (4.8), (4.9) and (4.12) we complete the proof.
\enddemo

\proclaim{Remark 4.1} In Proposition 3.1.1 of \cite{M1}, a
(considerably weaker) comparison on the relative volumes in the
similar spirit of Corollary 4.1 was established for the zero
divisors of holomorphic functions of polynomial growth, under
further assumptions on $M$ being of maximum volume growth and of
quadratic curvature decay, using very different method.
\endproclaim

The following consequence of Corollary 4.1 is somewhat surprising.
The result sharpen the Bishop-Gromov volume comparison theorem in
the presence of compact subvarieties.

\proclaim{Corollary 4.2} Let $M^m$ be a complete K\"ahler manifold
with nonnegative holomorphic bisectional curvature. Suppose that
$M$ contains a compact subvariety ${\Cal V}$ of complex dimension
$s$. Then there exists $C=C(m, s)>0$ such that  for $\delta(s)
\rho\ge \rho'\gg1$,
$$
\frac{V_{x_0}(\rho)}{V_{x_0}(\rho')}\le
C\left(\frac{\rho}{\rho'}\right)^{2(m-s)}.
$$
In particular,
$$
\lim_{\rho\to \infty}\frac{V_{x_0}(\rho)}{\rho^{2(m-s)}}<\infty.
$$
\endproclaim

An application of the above `monotonicity principle' we can have
another proof of Theorem 3.1 of \cite{N4}. In fact the new proof
gives a characterization of divisors defined by holomorphic
functions of polynomial growth (also called `polynomial functions'
according to the notation in \cite{W2}).

\proclaim{Theorem 4.1} Let $M$ be a complete K\"ahler manifold
with nonnegative bisectional curvature. Let ${\Cal V}$ be a
analytic divisor of $M$. Define the Lelong number (elliptic)  at
infinity of ${\Cal V}$ by
$$
\nu_{\infty} ({\Cal V})=\sup_{x_0\in M}\limsup_{\rho \to
\infty}\frac{\pi\rho^2 {\Cal A}_{{\Cal V},
x_0}(\rho)}{V_{x_0}(\rho)}. \tag 4.13
$$
Assume further that $H^1(M, {\Cal O}^*)=0$. Then
 ${\Cal V}$ is defined by a polynomial function if and only if
 $\nu_{\infty}({\Cal V})<\infty$.
Moreover, if ${\Cal V}=Z(f)$ for some $f\in P_d(M)$ (the space of
holomorphic functions of polynomial growth with degree at most
$d$) then
 there exists a $C(m)$ such that for any $x_0\in M$,
the Lelong number $\nu(x_0, {\Cal V})$ at $x_0$ is bounded by
$$
\nu(x_0, {\Cal V})\le C(m)\nu_{\infty} ({\Cal V})\tag 4.14
$$
and
$$
\nu_{\infty} ({\Cal V})\le C(m)\, d. \tag 4.15
$$
\endproclaim
\demo{Proof} Notice that (4.14) follows from Corollary 4.1
directly. First assume that ${\Cal V}$ is the zero divisor of a
polynomial function $f$ ($\in P_d(M)$ for some $d$). We shall show
that $\nu_{\infty}({\Cal V})<\infty$ (in fact (4.15)). We follow
the notation in Section 3 of \cite{N4} (also Section 2 of this
paper). Let $v(x,t)=\int_M H(x,y,t)\D \log |f|^2\, dv_y$ and
$w(x,t)=\frac{\p}{\p t}v(x, t)$. By Corollary 2.1 we know that
$$
tw(x,t)=(\pi t)\int_{\Cal V}H(x,y,t)\, dA_{\Cal V}.
$$
By Li-Yau's  lower bound estimate on heat kernel we have that
$$
tw(x,t)\ge C(m)\frac{\pi t {\Cal A}_{{\Cal V},
x}(\sqrt{t})}{V_x(\sqrt{t})}.
$$
Then (4.15) follows from (3.13) of \cite{N4}.

Now we assume that $\nu_\infty({\Cal V})<\infty$ we prove that
${\Cal V}$ is the divisor of a polynomial function. First, by the
solution to Cousin problem II (directly from the vanishing of the
 cohomology $H^1(M, {\Cal O}^*)$) we know that there exists a holomorphic function $f$ such that
$Z(f)={\Cal V}$. First we apply the `moment type estimate" from
\cite{N1} to estimate $tw(x,t)$ from above. Note that
$w(x,t)=\int_{\Cal V}H(x,y,t)\, dA_{\Cal V}(y)$ is well-defined
due to the assumption that $\nu_{\infty}({\Cal V})<\infty$. By
Corollary 2.1 we know that
$$
w(x,t)=\int_M H(x,y,t)\D \log |f|^2\, d\mu(y).
$$
Applying the `moment type estimate', Theorem 3.1 of \cite{N1} we
have that
$$
tw(x,t)\le C(m)\nu_{\infty}({\Cal V}).\tag 4.16
$$
Now we define $v(x,t)=\int_0^t w(x,\tau)\, d\tau+\log |f|^2(x)$.
Then (4.16) implies that $v(x,t)\le C(m)\nu_{\infty}({\Cal V})\log
(t+1) +\log |f|^2(x)$ for $t\gg 1$. One can also check that
$v(x,t)$ is a solution to the heat equation $\heat v(x,t)=0$ with
$v(x,0)=\log |f|^2(x)$. Now we apply the `moment type estimate' of
\cite{N1} again we have that
$$
\aint_{B_x(r)} \log |f|^2(y)\, d\mu(y) \le C(m)\nu_{\infty}({\Cal
V})\log (r+1).
$$
Now applying the mean-value inequality of Li-Schoen \cite{LS} we
conclude that $f$ is of polynomial growth.
\enddemo

\proclaim{Remark 4.2} In \cite{St} (see also \cite{Ru} for the
codimension one case), Stoll proved that an analytic divisor
${\Cal V}$ in $M=\C^m$ is algebraic if and only if
$\nu_{\infty}({\Cal V})<\infty$. The result was later generalized
to the case $M$ being an affine algebraic variety in \cite{GK} by
Griffiths and King. In fact,  in \cite{St} the result was proved
for any analytic sets of $\C^m$. It is desirable to generalize our
result to the high codimension case.

The assumption on vanishing of the cohomology in the above theorem
is satisfied, for example, when $M$ is Stein (cf. Theorem 5.5.2 of
\cite{H\"o}).
\endproclaim

In \cite{N4}, the author developed a parabolic approach to compare
the vanishing order of a holomorphic function, at any fixed point,
with the growth order at infinity. The method is sharp and
effective. It turns out that the parabolic method there is also
related to the Nevanlinna theory for several complex variables. In
order to illustrate this connection we need to recall some basic
notations from the Nevanlinna theory (cf. \cite{Gr}).

We call that a function $s: \R_{+}\to \R_{+}$ has finite order if
$$
Ord(s):=\limsup_{r\to \infty}\frac{\log s(r)}{\log r}<\infty.
$$
For a $f\in {\Cal O}(M)$, the space of holomorphic functions, we
define the order of $f$ in sense of Hadamard by
$$
Ord_H(f)=Ord(\log (A(r)))
$$
where $A(r)=\sup_{x\in B_o(r)}|f|(x)$ with $o\in M$ being a fixed
point. Following \cite{Gr, St} for any analytic subvariety of
complex dimension $s$ we define
$$
n_{\Cal V} (x_0, r)=\frac{{\Cal A}_{x_0}(r)(\pi
r^2)^{m-s}}{V_{x_0}(r)}
$$
and
$$
N_{\Cal V}(x_0, r)=\int_0^r \left(n_{\Cal V}(x_0, \tau)-n_{\Cal
V}(x_0, 0)\right)\frac{ d\tau}{\tau}+n_{\Cal V}(x_0, 0)\log r.
$$
$N_{\Cal V}(x_0, r)$ ($n_{\Cal V}(x_0, r)$) is called the counting
function since for $m=1$ and ${\Cal V}$ being the zeros of a
holomorphic function, it simply counts the number of zeros in
$B(x_0, r)$. If $x_0$ does not lie in ${\Cal V}$ one has $N_{\Cal
V}(x_0, r)=\int_0^r n_{\Cal V}(x_0, \tau)\frac{ d\tau}{\tau}$. One
can view the estimate (4.15) as bounding the counting function by
the growth order, a Nevanlinna type inequality. The following
result is a further generalization of (4.15), whose correspondence
in the Euclidean space is also known as the transcendental
B\'ezout estimate.

\proclaim{Corollary 4.3} Let $M$ be a complete noncompact K\"ahler
manifold with non-negative Ricci curvature. Let $f\in {\Cal O}(M)$
be a holomorphic function of finite order. Let $Z(f)$ be the zero
divisor of $f$. Then
$$
Ord(N_{Z}(x_0, r))\le Ord_H(f). \tag 4.17
$$
\endproclaim
\demo{Proof} The proof follows the same line of argument as the
proof of Theorem 4.1. We leave it to the interested readers.
\enddemo
We found the connection between the parabolic equations (as well
as the related differential Harnack inequality) and the Nevanlinna
theory quite interesting.

\input amstex
\documentstyle{amsppt}
\magnification=1200 \hsize=13.8cm \catcode`\@=11
\def\NoLogo{\let\logo@\empty}
\catcode`\@=\active \NoLogo

\def\heat{\lf(\frac{\p}{\p t}-\Delta\ri)}

\def \b {\beta}

\def\Ric{\text{Ric}}
\def\lf{\left}
\def\ri{\right}
\def\bbar{\bar \beta}
\def\a{\alpha}

\def\g{\gamma}
\def\e{\epsilon}
\def\p{\partial}
\def\delbar{{\bar\delta}}

\def\dbar{\bar\partial}

\def\ctL{{\Cal L}_+}
\def\C{\Bbb C}
\def\R{\Bbb R}
\def\P{\Bbb P}
\def\E{\Bbb E}

\def\vp{\varphi}
\def\test{\vp_{(x_0, t_0), \rho}}
\def\testt{\vp_{(x_0, t_1), \rho}}

\def\tN{\tilde N}
\def\cN{\Cal N}
\def\tZ{\tilde Z}
\def\ctN{\tilde{\Cal N}}

\def\dbar{\bar\partial}

\def\bb{{\bar\beta}}
\def\abb{{\alpha\bar\beta}}
\def\gbd{{\gamma\bar\delta}}

\def \D {\Delta}
\def\aint{\frac{\ \ }{\ \ }{\hskip -0.4cm}\int}
\documentstyle{amsppt}
\vsize=19.0 cm

\subheading{\S5 Heat kernel and the reduced volumes}

\vskip .2cm

In this section we consider the fundamental solution to the
time-dependent heat equation with metrics deformed by
K\"ahler-Ricci/Ricci flow (1.1). As we explained in Section 2 that
the fundamental solution $H(x,y,t, t_0)$ satisfies (1.2). We shall
prove a sharp lower bound on $H(x,y,t,t_0)$. In \cite{N5} we
derived an estimate for the fundamental solution satisfying the
time-dependent heat equation itself (instead of (1.2)). The
estimate there is only valid for short time interval. In
\cite{Gu}, a rough lower bound was obtained through the earlier
mentioned Harnack inequality, Theorem 1.3. Notice that the result
in \cite{Gu} is not  sharp and the result in \cite{N5} is only
sharp in the exponents. In this section we show a  sharp lower
bound on the heat kernel in the case that $M$ is  either a
complete Riemannian bounded nonnegative curvature operator or a
complete K\"ahler manifold with bounded nonnegative bisectional
curvature. For the sake of simplicity we only state the result for
the K\"ahler-Ricci flow and leave the Riemannian analogue to the
interested readers. Before we state our result we need to recall
some notations and computations from \cite{FIN} (which follows
closely  the computation in \cite{P}). For the simplicity we
assume that $t_0=0$.

Let $g(t)$ solve the K\"ahler-Ricci flow on $M^m\times [0, T]$
($m=\dim_{\C}(M)$ and $n=2m$). Fix $x_0$ and let $\gamma$ be a
path $(x(\eta),\eta)$ joining $(x_0, 0)$ to $(y, t)$. Following
\cite{P} and \cite{LY, FIN} we define
$$
\ctL(\gamma)=\int_0^t\sqrt{\eta} \left({\Cal
R}+4|\gamma'(\eta)|^2\right)\, d\eta.\tag 5.1
$$
Let $X=\gamma'(t)=\frac{d z^{\a}(t)}{d t}\frac{\p}{\p z^{\a}}$ and
let $Y$ be a variational vector field along $\gamma$. Here
$|\gamma'(t)|^2=g_{\abb} \frac{d z^{\a}(t)}{d t}\frac{d
z^{\bar{\beta}}(t)}{d t}$. Using $\ctL$ as energy we can define
the $\ctL$-geodesics and we denote $L_{+}(y, t)$ to be the length
of a shortest geodesics jointing $(x_0, 0)$ to $(y, t)$. We also
define
$$
\ell_+(y, t):=\frac{1}{2\sqrt{t}}L_+(y, t).
$$
Following the first and second variation calculation of \cite{P}
(see also \cite{FIN}) we compute that
$$
\gather |\nabla \ell_+|^2=-{\Cal
R}+\frac{\ell_+}{t}+\frac{K}{t^{3/2}}, \tag
5.2\\
\frac{\p \ell_+}{\p t} ={\Cal R}
-\frac{K}{2t^{3/2}}-\frac{\ell_+}{t},
\tag 5.3\\
\D \ell_+\le {\Cal R}+\frac{n}{2t} -\frac{K}{2t^{3/2}}. \tag 5.4
\endgather
$$
Here
$$
K:=\int_0^t \eta^{3/2}H(X)\, d\eta,
$$
where $H(X):=\p{\Cal R}/\p t+2\langle\nabla {\Cal R},
X\rangle+2\langle X, \nabla {\Cal R}\rangle +4\Ric(X,X)+{\Cal
R}/t$ is  exactly  the traced Li-Yau-Hamilton differential Harnack
expression in \cite{C1} applying to the $(1,0)$ vector field $2X$.

\proclaim{Theorem 5.1} Let $(M^m, g(t))$ be a complete solution to
K\"ahler-Ricci (Ricci) flow with bounded nonnegative bisectional
curvature (curvature operator). Let $H(x_0, y, t)$ be the
fundamental solution to the time dependent heat equation
(satisfying (1.2) and $H(x_0, y, 0)=\delta_{x_0}(y)$). Then
$$
\tilde u(x,t):= \frac{1}{(\pi t)^{m}}\exp\left(-\ell_+(x,
t)\right)
$$
satisfies
$$
\left(\frac{\p}{\p t}-\Delta -{\Cal R}\right)\tilde u (x,t)\le 0.
\tag 5.5
$$
In particular,
$$
\tilde u(x,t)\le H(x_0, y, t) \tag 5.6
$$
and
$$
\tilde \theta_{+}^{(x_0, 0)}(t):=\int_M \tilde u(x,t)\, d\mu_t
$$
is monotone decreasing.
 Moreover,  the equality in (5.5), or (5.6) implies that $M$ is a
gradient expanding soliton.
\endproclaim
\demo{Proof}  First (5.2)--(5.4) implies that $$\left(\frac{\p}{\p
t}-\Delta -{\Cal R}\right) \left(\frac{1}{(\pi
t)^{m}}\exp\left(-\ell_+(y,
t)\right)\right)=-\frac{K}{t^{\frac{3}{2}}}\left(\frac{1}{(\pi
t)^{m}}\exp\left(-\ell_+(y, t)\right)\right)\le 0.
$$
Here we have used fact that $K\ge 0$ under the assumption that $M$
has bounded non-negative bisectional curvature. Also if the
equality holds it implies that $K\equiv 0$. This further implies
that $M$ is an expanding soliton from the computation in
\cite{FIN}. In order to prove (5.6) one just need to apply the
maximum principle (cf. \cite{NT2, N5}) and  notice that
$\lim_{t\to 0}\frac{1}{(\pi t)^m}\exp(-\ell_+(y,
t))=\delta_{x_0}(y)$. The equality case follows from the
consideration of the equality in \cite{FIN}.
\enddemo

\proclaim{Remark 5.1} In \cite{CY}, Cheeger and Yau proved that
for complete manifolds with nonnegative Ricci curvature, the heat
kernel $H(x,y, t)$ has the lower bound estimate:
$$
\bar{H}(d(x, y), t):=\frac{1}{(4\pi
t)^{\frac{n}{2}}}\exp(-\frac{d^2(x,y)}{4t})\le H(x,y, t)
$$
by showing that the transplant of the Euclidean heat kernel $\bar
H(d(x, y), t)$ is a sub-solution to the heat equation. Recall that
in \cite{P}, Perelman first discovered that there is a similar
result for the backward Ricci flow, even without any {\it
curvature sign} assumptions (this is the astonishing part of
Perelman's work). Namely he proved, for $\frac{\p}{\p \tau}
g_{ij}=2R_{ij}$ and the reduced distance $\ell (y, \tau)$, defined
formally by the same expression as $\ell_+$, that
$$
\left(\frac{\p}{\p \tau} -\D +{\Cal R}\right)\bar u(y, \tau) \le 0
\tag 5.7
$$
where
$$
\bar u(y, \tau)=\frac{1}{(4\pi
\tau)^{\frac{n}{2}}}\exp(-\ell(y,\tau)).
$$
In particular, Perelman further showed (in Corollary 9.5 of
\cite{P})  that $\bar u $ gives a lower bound for the heat kernel
of the {\bf backward heat equation} ($(\frac{\p}{\p \tau}-\D)\bar
v=0$). (Notice that the heat kernel $H(x,y, \tau)$ for the
backward heat equation satisfies the Schr\"odinger eqaution
$\left(\frac{\p}{\p \tau} -\D +{\Cal R}\right)H(x,y, \tau)=0$. One
should also refer to \cite{LY, Theorem 4.3} for a precedence of
Perelman's second variation computation in \cite{P, Section 7}.)
Thus  both Perelman's monotonicity of the reduce volume and
Theorem 5.1 above can be viewed as nonlinear analogue of the
earlier work of Cheeger-Yau in \cite{CY} and Li-Yau in \cite{LY}.
\endproclaim

 Tracing the computation of \cite{FIN} we also have the following
estimates for $\tilde u(x,t)$.

\proclaim{Proposition 5.1} Assume that $(M, g(t))$ be a
K\"ahler-Ricci flow with bounded nonnegative bisectional curvature
on $M\times [0, T)$. Let $\tilde u(x,t)$ be as above. Then
$$
\log (\tilde u)_{\abb} +R_{\abb}+\frac{1}{t}g_{\abb}\ge 0. \tag
5.8
$$
The equality holds if and only if $(M, g(t))$ is an expanding
K\"ahler-Ricci solition.
\endproclaim
Notice that (5.7) is equivalent to  the estimate (1.3) on $u(x,t)$
a  positive solution to the {\it forward conjugate heat equation}.
On the other hand $\tilde u(x,t)$ is only a sub-solution to the
{\it forward conjugate heat equation}.

\proclaim{Remark 5.2} In \cite{FIN, Corollary 2.1} we showed that
under the assumption of bounded nonnegative curvature operator
there exists a same estimate as {\rm (5.8)} (without the bars).
This suggests that one may have the similar estimate as Theorem
1.2 for the Riemannian case on positive solution $u$ of the {\it
forward conjugate heat equation}. Indeed, this can be shown if one
proves,  under the assumption of $M$ having  bounded nonnegative
curvature operator, that
$$
Y_{ij}:= \nabla_i\nabla_j {\Cal R}-2\nabla_k R_{ij}\nabla_k(\log
u)+2R_{ikjl}\nabla_k(\log u)\nabla_l(\log u)
 +R_{ik}R_{jk}+\frac{1}{t}R_{ij}\ge 0. \tag 5.9
$$
This claim is based on the  following computation which is
essentially due to Chow. Note that $Y_{ij}\equiv 0$ on gradient
expanding solitons and its trace is the same as the trace of
Hamilton's matrix Li-Yau-Hamilton expression in \cite{H1}.
\endproclaim

\proclaim{Lemma 5.1}  Let
$$\tZ_{ij}=R_{ij}+(\log
u)_{ij}+\frac{1}{2t}g_{ij}.$$ Then
$$
\split \left(\frac{\p}{\p t}-\D_L \right)
\tZ_{ij}&=Y_{ij}+2\nabla_k \tZ_{ij}\nabla_k (\log
u)\\
&\quad
+\tZ_{ik}\left(\log(u)_{jk}-R_{jk}-\frac{1}{2t}g_{jk}\right)+
\left(\log(u)_{ik}-R_{ik} -\frac{1}{2t}g_{ik}\right)\tZ_{jk}
\endsplit \tag 5.10
$$
where $\D_L\eta_{ij}=\D
\eta_{ij}+2R_{ikjl}\eta_{kl}-R_{ik}\eta_{kj}-R_{ik}\eta_{ik}$ is
the Lichnerowicz operator acting on the symmetric $2$-tensor
$\eta_{ij}$.
\endproclaim

The above computation was first carried out in \cite{Ch2} for the
{\it backward Ricci flow}, with $Y_{ij}$ is replaced by Hamilton's
matrix Harnack expression for shrinkers. Notice also that the
matrix Harnack expression $Y_{ij}$ and Hamilton's expression for
Ricci expanders are the same if the manifold is K\"ahler.

One can obtain some upper bound on the heat kernel $H(x,y,t, 0)$
using the Harnack inequality proved in Corollary 2.1 and the fact
that $\int_M h(x,y,t, 0)\, d\mu_t=1$. But the result is not as
satisfactory as for the fixed metric case of Li-Yau. Hence we
shall leave this to a later investigation. On the other hand, the
localization technique of Ecker can be applied to the reduced
volume for the Ricci expanders (defined in \cite{FIN}) to obtain
the  monotonicity of a localized reduced volume {\it without
assuming bisectional curvature (nor curvature operator) being
nonnegative}. Recall that in \cite{FIN}, the authors proved that,
if $M$ is a closed manifold, the {\it forward reduced volume}
$$ \theta_+^{(x_0, 0)}(t):=\int_M
\hat u(x,t)\,d\mu_t
$$
is monotone non-increasing, where $ \hat u(x,t)=\frac{e^{\ell_+(x,
t)}}{(\pi t)^m}. $ Here we define $\hat{u}$  for the K\"ahler case
to be coherent with our previous discussions. This monotonicity
has severe restriction since the reduced volume $\theta_{+}^{(x_0,
0)}(t)$ is only meaningful when $M$ is compact. We shall show the
monotonicity of a localized reduced volume, which is well-defined
on complete manifolds, to compensate such restriction. In order to
put the result in its general form we denote by $l_+^{(x_0, t_0)}$
the {\it forward reduced distance} with respect to $(x_0, t_0)$ if
one replaces the reference point $(x_0, 0)$ by $(x_0, t_0)$.
Correspondingly we denote by $\hat u^{(x_0, t_0)}$ and
$\theta_{+}^{(x_0, t_0)}$, the reduced volume function and the
reduced volume  with respect to $(x_0, t_0)$, respectively. Let
$$
\bar{L}^{(x_1, t_1)}_+=(t-t_1)\ell_+^{(x_1, t_1)}
$$
and
$$
\vp^{(x_1, t_1)}_{t_2, \rho}(x, t)=\left(1-\frac{\bar{L}_+^{(x_1,
t_1)}(x,t)+m(t-t_2)}{\rho^2}\right)_{+}.
$$
 Then (5.2)--(5.4) imply that
$$
\lf(\frac{\p}{\p t}+\D -{\Cal R}\ri) \hat u^{(x_0, t_0)}(x,t)\le 0
\tag 5.11
$$
and
$$
\heat \vp^{(x_1, t_1)}_{t_2, \rho}(x, t)\le 0.\tag 5.12
$$
Notice that if ${\Cal R}$ is uniformly bounded from below,
$\vp^{(x_1, t_1)}_{t_2, \rho}(x, t)$ is compactly supported. We
then have the following localized monotonicity of the reduced
volume. \proclaim{Proposition 5.2}
$$
\frac{d}{d t}\theta^{(x_0, t_0)}_{+, \vp}(t) \le 0\tag 5.13
$$
where $$ \theta^{(x_0, t_0)}_{+, \vp}(t)=\int_M\vp^{(x_1,
t_1)}_{t_2, \rho}(x, t)\hat u^{(x_0, t_0)}(x,t)\, d\mu_t. \tag
5.14
$$
\endproclaim

It turns out that one can also construct a subsolution to the heat
equation with compact support, similar to $\vp^{(x_1, t_1)}_{t_2,
\rho}(x, t)$ defined above. This is through Perelman's {\it
reduced distance} $\ell^{(x_1, t_1)}(x, \tau)$, with $\tau
=t_1-t$. ($\ell^{(x_1, t_1)}(x, t_1-t)$ is only defined for $t\le
t_1$. Please refer to \cite{P, Section 7} for the detailed
discussions.) Here we use $(x_1, t_1)$ to specify the reference
space-time point with respect to which the reduced distance is
defined. Recall that from \cite{P, Section 7}, $\bar{L}^{(x_1,
t_1)}(x, t):=(t_1-t)\ell^{(x_1, t_1)}(x, t_1-t)$ (again we do not
have the factor $4$ since we are in the K\"ahler setting)
satisfies the differential inequality:
$$
\left(-\frac{\p}{\p t}+\D\right)\bar{L}^{(x_1, t_1)}(x, t)\le m.
$$
Then in the case ${\Cal R}$ is bounded from below, we may define a
compact supported function
$$
\psi^{(x_1, t_1)}_{t_2, \rho}(x, t):=\left(1-\frac{\bar{L}^{(x_1,
t_1)}(x,t)+m(t-t_2)}{\rho^2}\right)_{+}.
$$
It is easy to see that
$$
\heat \psi^{(x_1, t_1)}_{t_2, \rho}(x, t)\le 0.\tag 5.15
$$
This gives another localization on the {\it forwarded reduced
volume}.

\proclaim{Proposition 5.3}$$ \frac{d}{d t}\theta^{(x_0, t_0)}_{+,
\psi}(t) \le 0\tag 5.16
$$
where $$ \theta^{(x_0, t_0)}_{+, \vp}(t)=\int_M\psi^{(x_1,
t_1)}_{t_2, \rho}(x, t)\hat u^{(x_0, t_0)}(x,t)\, d\mu_t. \tag
5.17
$$
\endproclaim

The similar idea can also be applied to Perelman's entropy. Let
$(x_0, t_0)$ be a fixed space-time point. Let $\tau=t_0-t$ and
$u(x, \tau)$ be the fundamental solution of the {\it backward
conjugate heat equation} $\frac{\p}{\p \tau}-\D +{\Cal R}$
centered at $(x_0, t_0)$. Write $u=\frac{e^{-f}}{(\pi \tau)^{m}}$
and define $v=\left(\tau\left(2\D f-|\nabla f|^2+{\Cal
R}\right)+f-2m\right)u$. Then it was proved in \cite{P} that
$$
\left(\frac{\p}{\p \tau}-\D +{\Cal R}\right)v=-\tau
\left(|R_{\abb}+\nabla_{\a}\nabla_{\bar{\b}}f-\frac{1}{\tau}g_{\abb}|^2+|\nabla_{\a}\nabla_{\b}f|^2\right)u
\tag 5.18
$$
and $v\le 0$. Applying (5.12) or (5.15) we can have the following
local monotonicity formulae.

\proclaim{Proposition 5.4}
$$
\frac{d}{d t}\left(\int_M -v\vp^{(x_1, t_1)}_{t_2, \rho}\,
d\mu_t\right)\le -\int_M \tau
\left(|R_{\abb}+\nabla_{\a}\nabla_{\bar{\b}}f-\frac{1}{\tau}g_{\abb}|^2+|\nabla_{\a}\nabla_{\b}f|^2\right)
u\vp^{(x_1, t_1)}_{t_2, \rho}\, d\mu_t \tag 5.19
$$
and
$$
\frac{d}{d t}\left(\int_M -v\psi^{(x_1, t_1)}_{t_2, \rho}\,
d\mu_t\right)\le -\int_M \tau
\left(|R_{\abb}+\nabla_{\a}\nabla_{\bar{\b}}f-\frac{1}{\tau}g_{\abb}|^2+|\nabla_{\a}\nabla_{\b}f|^2\right)
u\psi^{(x_1, t_1)}_{t_2, \rho}\, d\mu_t . \tag 5.20
$$
\endproclaim
\demo{Proof} The proof is just direct computations and integration
by parts, using (5.12), (5.15) and (5.18).
\enddemo

There is no difference for the Riemannian cases. The formulae are
exactly the same except some factors caused by the definition of
the operators.

We conclude this section by an application of Theorem 5.1 to the
study of the large time behavior of K\"ahler-Ricci flow. Let
$(M^m, g_0)$ be a complete K\"ahler manifold with bounded
nonnegative holomorphic bisectional curvature. If we further
assume that $(M, g_0)$ has maximum volume growth, namely for any
$x$, the volume of the ball $B(x, r)$, $V_x(r)$ is bounded from
below by $\delta r^{2m}$ for some positive $\delta$, in \cite{N6}
we proved that the K\"ahler-Ricci flow (1.1) has long time
solution with the initial data $g(x, 0)=g_0(x)$. Moreover, there
exists a constant $A=A(M)>0$ such that  $t\gg 1$,
$$
t{\Cal R}(x, t)\le A. \tag 5.21
$$
Applying Theorem 5.1 to this situation we can show the following
result.

\proclaim{Corollary 5.1} Let $(M, g_0)$ be a complete K\"ahler
manifold as above. Let $g(x, t)$ be a long time solution to
K\"ahler-Ricci flow as above. For any $(x_j, t_j)$ with $t_j \to
\infty$, define $g_j (t)=\frac{1}{t_j}g(t_j t)$. Then the pointed
sequence $(M, x_j, g_j(x, t))$ sub-sequentially converges to a
gradient expanding K\"ahler-Ricci soliton $(M_\infty, x_\infty,
g_\infty(t))$.
\endproclaim

\demo{Proof} For any $t$, denote by $B_t(x, r)$ the ball of radius
$r$ with respect to $g(t)$, and by $V_t(x, r)$ the volume of this
ball (with respect to $g(t)$). By Theorem 2.2 of \cite{NT3}, we
know that $V_t(x, r)\ge \delta r^{2m}$. Together with (5.21) we
have the injectivity radius bounded from below uniformly for
$g_j(t)$. Therefore, by Hamilton's  compactness theorem, $(M, x_j,
g_j(t))$ subsequentially converges to $(M_\infty, x_\infty,
g_\infty(t))$, a solution to K\"ahler-Ricci flow defined on
$M_\infty \times (0, \infty)$. The only thing we need to prove is
that $(M_\infty, g_\infty(t))$ is an expanding soliton. It is easy
to see that $(M_\infty, g_\infty)$ also has bounded nonnegative
bisectional curvature and the maximum volume growth. By Corollary
1 of \cite{N6}, we know that $M_\infty$ is topologically
$\R^{2m}$. The fact that it is an expanding soliton follows from
the monotonicity consideration. Let us adapt the notations from
Theorem 5.1. Let $(x_0, 0)$ be a fixed reference point, with
respect to which we have the {\it forward reduced distance
function} $\ell_+(x, t)$, and the {\it second forward reduced
volume  function} $\tilde u(x, t)$ and the {\it second forward
reduced volume} $\tilde \theta_+(t)$. By the proof of Theorem 5.1
we have that
$$
\frac{d }{d t} \tilde \theta_+(t)=-\int_M \frac{K}{t^{3/2}}\tilde
u\, d\mu_t,\, dt
$$
from which we can conclude that
$$
\tilde \theta_+(t_j)-\tilde \theta_+(2t_i)=\int_{t_j}^{2t_j}
\frac{1}{t^{3/2}}\int_M K \tilde u\, d\mu_t.
$$
Hence
$$
\tilde \theta_+(2t_j)-\tilde \theta_+(t_j) \ge \int_{t_j}^{2t_j}
\frac{1}{t^{3/2}}\int_M  K \tilde u\, d\mu_t \, dt.
$$
Taking limit we have that for the limit metric
$$
\int_{1}^{2} \frac{1}{t^{3/2}}\int_{M_\infty}  K \tilde u\, d\mu_t
\, dt=0.
$$
The conclusion now follows from \cite{C2} or Theorem 4.1 of
\cite{N2} since $K=0$ implies that the linear trace LYH quantity
$H$ achieves its minimum (zero) somewhere.
\enddemo

\proclaim{Remark 5.3} Corollary 5.1  simply says that the
blow-down limit of a Type III solution is an expanding soliton.
The result is in some sense dual to Proposition 11.2 of \cite{P}.
It is also similar to the situation studied in \cite{Hu} for the
mean curvature flow. The maximum volume growth assumption can be
weaken to a certain  $\kappa$-noncollapsing condition  in terms of
the lower bound of  $\tilde \theta_+(t)$.

It has been proved in \cite{CT1} (see also \cite{CT2} for further
more recent  progresses) that a gradient expanding K\"ahler
soliton must be biholomorphic to $\C^m$.

One can have a similar result for the Riemannian case, replacing
the nonnegativity of the bisectional curvature by the
nonnegativity of the curvature operator.

\endproclaim

\subheading{\S  Appendix: A parabolic relative volume comparison
theorem}\vskip .2cm

Recall that Cheeger and Yau proved that on a complete Riemannian
manifold with nonnegative Ricci curvature, the heat kernel $H(x,y,
\tau)$ (the fundamental solution of the operator  $\heatau $) has
the lower estimate
$$
H(x,y, \tau)\ge \frac{1}{(4\pi
\tau)^{\frac{n}{2}}}\exp\left(-\frac{r^2(x,y)}{4\tau}\right) \tag
A1
$$
where $r(x,y)$ is the distant function on the manifold. This fact
can be derived out of the maximum principle and the differential
inequality
$$
\heatau \left(\frac{1}{(4\pi
\tau)^{\frac{n}{2}}}\exp\left(-\frac{r^2(x,y)}{4\tau}\right)\right)
\le 0.
$$
Integrating on the manifold $M$, this differential inequality also
implies the monotonicity (monotone non-increasing) of the integral
$$
\tilde V(x_0, \tau):=\int_M\frac{1}{(4\pi
\tau)^{\frac{n}{2}}}\exp\left(-\frac{r^2(x_0,y)}{4\tau}\right)
\tag A2
$$
 In \cite{P}, Perelman discovered a striking analogue
of this comparison result for  the Ricci flow geometry. More
precisely, he introduced a length function, called the {\it
reduced distance},
$$\ell_{x_0, g(\tau)}( y, \taub)
:=\inf_{\gamma}\frac{1}{2\sqrt{\tau}} \int_0^{\taub}
\sqrt{\tau}|\gamma'(\tau)|^2\, d\tau$$ for all $\gamma(\tau)$ with
$\gamma(0)=x_0,\, \gamma(\taub)=y$ (in the right context, we often
omit the subscript $x_0, g(\tau)$),  and  a functional, called
{\it the reduced volume},
$$
\tilde V_{g(\tau)}(x_0, \tau):=\int_M\frac{1}{(4\pi
\tau)^{\frac{n}{2}}}\exp\left(-\ell(y, \tau)\right) \tag A2'
$$
with respect to a solution $g_{ij}(x, \tau)$ to the  backward
Ricci flow
$$
\frac{\p}{\p \tau} g_{ij}=2 R_{ij}
$$
on $M\times [0, a]$.   Perelman proved further that $\tilde
V_{g(\tau)}(x_0, \tau)$ is monotone non-increasing in $\tau$. In
fact, he
 showed this monotonicity result via a space-time relative
 comparison result and used it
  to give
a more flexible proof of the {\it  $\kappa$-noncollapsing} on  the
solution to Ricci flow, defined on $M\times[0, T]$,  for any given
finite time $T$ with given initial data.

In \cite{P}, Perelman also discovered a more rigid monotone
quantity, {\it the entropy functional}. These two functionals can
be related through a differential inequality of Li-Yau-Hamilton
type. One can refer to \cite{CLN} for an exposition on
 this relation. The entropy is
 stronger but less flexible.
The entropy monotonicity has its analogue for the positive
solutions to linear heat equation on a fixed Riemannian manifold
with non-negative Ricci curvature. This was derived  in \cite{N3},
where the author also observed that the entropy defined on a
complete Riemannian manifold $M$ with non-negative Ricci curvature
is  related to the volume growth of the manifold. More precisely,
let

$$
{\Cal W}(f, \tau):=\int_M \left(\tau |\nabla
f|^2+f-n\right)\frac{e^{-f}}{(4\pi \tau)^{\frac{n}{2}}}\, d\mu
\tag A3
$$
for any $\tau>0$ and $C^1$ function $f$ with
$\int_M\frac{e^{-f}}{(4\pi \tau)^{\frac{n}{2}}}\, d\mu=1$. If
$u(x, \tau)=\frac{e^{-f}}{(4\pi \tau)^{\frac{n}{2}}}$ is a
solution to the heat equation, it was proved in \cite{N3} that
${\Cal W}(f, \tau)$ is monotone non-increasing in $\tau$.
Moreover, it was shown that
$$
\lim_{\tau \to \infty}{\Cal W}(f, \tau)=\log \left(\lim_{\tau\to
\infty}\tilde V_{x_0}(\tau)\right)
$$
if $u(x, \tau)$ is a fundamental solution originated at $x_0$. In
fact, $$\lim_{\tau\to \infty}\tilde V_{x_0}(\tau)=\lim_{r\to
\infty}\frac{V_{x_0}(r)}{\omega_n r^n}$$ where $V_{x_0}(r)$ is the
volume of ball of radius $r$, which is independent of $x_0$. (The
limit on the right hand side above is also called the cone angle
at infinity).

Motivated by this close connection between the Ricci flow geometry
for a family of metrics  and the Riemannian geometry of a fixed
Riemannian metric, it is nature to seek a localized version of the
above mentioned Cheeger-Yau's result on the monotonicity of the
reduced volume $\tilde V_{x_0}(\tau)$. In fact,  a local version
of the heat kernel comparison was  also carried out in the paper
of \cite{CY} from a PDE point of view. However,  Perelman's
localization  in the case of Ricci flow geometry is along the line
of comparison geometry, and very much different from Cheeger-Yau's
localization consideration (by considering the Dirichlet or
Neumann boundary value problem for the heat equation). This leads
us to formulate a  new  relative (local) volume comparison theorem
in this section, which is the linear analogue of Perelman's
formulation for Ricci flow. It can be viewed a parabolic version
of the (classical) relative volume comparison theorem in the
standard Riemannian geometry. The interested readers may want to
compare these two results in more details. The relation between
the classical relative volume comparison theorem and the new one
here is very similar to the one between the monotonicity formula
for the minimal submanifolds in $\R^n$ and Huisken's monotonicity
for mean curvature flow in $\R^n$ (as well as Ecker's localized
version).

In order to state our result, let us first fix some notations. It
is our hope that the exposition below is detailed enough to be
helpful in  understanding \cite{P} better.  Fix a  point $x_0\in
M$. Let $\gamma(\tau)$ ($0\le \tau\le \bar{\tau}$) be a curve
parameterized by the time variable $\tau$ with $\gamma(0)=x_0$.
Here we image that we have a time function $\tau$, with which some
parabolic equation is associated. Define the ${\Cal L}$-length by
$$ {\Cal L}(\gamma)(\taub)=\int_0^{\taub}
\sqrt{\tau}|\gamma'(\tau)|^2\, d\tau. \tag A4
$$
We can define the ${\Cal L}$-geodesic to be the curve which is the
critical point of ${\Cal L}(\gamma)$. The simple computation shows
that the first variation of ${\Cal L}$ is given by
$$
\delta {\Cal L}(\gamma)=2\sqrt{\taub}\langle Y, X\rangle
(\taub)-2\int_0^{\taub}\sqrt{\tau}\left(\langle \nabla_X X
+\frac{1}{2\tau}X, Y\rangle\right)\, d\tau, \tag A5
$$
where $Y$ is the variational vector field, from which one can
write down the ${\Cal L}$-geodesic equation. It is an easy matter
to see that $\gamma$ is a ${\Cal L}$-geodesic if and only if
$\gamma(\sigma)$ with $\sigma =2\sqrt{\tau}$ is a geodesic. In
another word, a ${\Cal L}$-geodesic is a geodesic after certain
re-parametrization. Here we insist all curves are parameterized by
the `time'-variable $\tau$. One can check that for any $v\in
T_{x_0}M$ there exists a ${\Cal L}$-geodesic $\gamma(\tau)$ such
$\frac{ d}{d \sigma}\left(\gamma(\sigma)\right)|_{\sigma =0}=v.$
Notice that the variable $\sigma$  scales in the same manner as
the distance function on $M$. So it is more convenient to work
with $\sigma$.

We then define the ${\Cal L}$-exponential map by
$$
{\Cal L}\exp_{v}(\sigmab):=\gamma_{v}(\sigmab)
$$
if $\gamma_v(\sigma)$ is a ${\Cal L}$-geodesic satisfying that
$$
\lim_{\sigma \to 0}\frac{d
}{d\sigma}\left(\gamma_v(\sigma)\right)=v.
$$
It is also  illuminating to go one dimensional higher by
considering the manifold $\tM=M\times [0, 2\sqrt{T}]$ and the
space-time exponential map $\texp(\tilde v^a)=({\Cal
L}\exp_{\frac{v^a}{a}}(a), a)$, where $\tilde v^a=(v^a, a)$.
Denote  $\frac{v^a}{a}$ simply by $v^1$ and $(v^1, 1)$ by $\tilde
v^1$.  Also let $\tilde \gamma_{\tilde v^a}(\eta)=\texp(\eta\tilde
v^a)$. It is easy to see that
$$
\tilde \gamma_{\tilde v^a}(\eta)=\tilde \gamma_{\tilde v^1}(\eta
a).
$$
This shows that
$$
d\texp|_{(0, 0)}=\text{identity}.
$$
Computing the second variation of ${\Cal L}(\gamma)$ gives that
$$
\delta^2 {\Cal L}(\gamma)=2\sqrt{\taub}\langle \nabla_Y Y,
X\rangle+\int_0^{\taub}2\sqrt{\tau}\left(|\nabla_X Y|^2-R(X, Y, X,
Y)\right)\, d\tau \tag A6
$$
where $Y$ is a given variational vector field. Let $(y,
\sigmab)=({\Cal L}\exp_{v}(\sigmab), \sigmab)$. Consider the
variation which is generated by $(w,0)$. Namely consider the
family
$$
U(s, \sigma)=({\Cal L}\exp_{v+sw}(\sigma), \sigma)=\texp((\sigma
(v+sw), \sigma)).
$$
Direct calculation shows that
$$
\frac{D U}{\partial \sigma}(0, 0)=(v,1), \quad \frac{D U}{\partial
s}(0, 0)=0
$$
and that the Jacobi field $\tilde {J}_w(\sigma)=(J_w(\sigma), 0)$
is given by
$$
\frac{D U}{\p s}|_{s=0}(\sigma)=d\texp ((\sigma w, 0))
$$
with the initial velocity
$$
\nabla_{\frac{\p}{\p \sigma}}\left(\frac{D U}{\p
s}|_{s=0}\right)(0)=(w, 0).
$$
(One can define the Jacobi operator to be the linear second order
operator associated with the quadratic form in the right-hand side
of (A6). One call a vector field along $\gamma$ a ${\Cal
L}$-Jacobi field if it satisfies the Jacobi equation. It is  easy
 to show that the variational vector field of a family of ${\Cal
L}$-geodesics satisfies the Jacobi equation as in the standard
Riemannian geometry. In fact, the ${\Cal L}$-Jacobi field turns
out to be just the regular Jacobi-field after re-parametrization.)
This shows that
$$
(d {\Cal L}\exp)_{w}(\sigma)= J_w(\sigma) \tag A7
$$
and
$$
d\texp ((\sigma w, 0))=\tilde J_w(\sigma). \tag A8
$$
Now we can conclude that $(y, \sigma)$ is a regular value  of the
map $\texp$ (and $y$ is a regular value of ${\Cal
L}\exp_{(\cdot)}(\sigma)$) if and only if that any Jacobi field
$\tilde J$ with initial condition as above does not vanish at
$\sigma$. We can introduce the concept of conjugate point (with
respect to $x_0$) similarly as in the classical case. We can
define the set
$$
D(\sigmab)\subseteq T_{x_0}M
$$to be  the collection of vectors $v$ such that $({\Cal L}\exp_v(\sigma),
\sigma)$ is a ${\Cal L}$-geodesic along which there is no
conjugate point up to $\sigmab$. Similarly we can define the set
$$
C(\sigmab)\subseteq T_{x_0}M
$$
to be the collection of vectors $v$ such that $({\Cal
L}_v(\sigma), \sigma)$ is a minimizing ${\Cal L}$-geodesic up to
$\sigmab$. One can see easily that $D(\sigma)$ and $C(\sigma)$
decreases (as sets) as $\sigma$ increases. For any measurable
subset $A\subseteq T_{x_0}M$ we can define
$$
D_A(\sigma)=A\cap D(\sigma) \quad \text{and}\quad
C_A(\sigma)=A\cap C(\sigma).
$$
Following \cite{P} we also introduce the $\ell$-`distance'
function.
$$
\ell_{x_0}(y, \taub)=\frac{1}{2\sqrt{\taub}}L_{x_0}(y, \taub),
\quad \text{where }L_{x_0}(y, \taub)=\inf_{\gamma}{\Cal
L}(\gamma)).
$$
Here our $\ell$ is defined for a fixed background metric. We also
omit the subscript $x_0$ in the context where the meaning is
clear. The very same consideration as in \cite{P}, as well as in
the standard Riemannian geometry, shows that
$$
|\nabla \ell|^2=\frac{1}{\taub}\ell\tag A9
$$
$$
\ell_\tau=-\frac{1}{\taub}\ell \tag A10
$$
and
$$
\Delta \ell \le
\frac{n}{2\taub}-\frac{1}{\taub^{\frac{3}{2}}}\int_0^{\taub}
\tau^{\frac{3}{2}}\Ric(X, X)\, d\tau \tag A11
$$
where $X=\gamma'(\tau)$ with $\gamma(\tau), \, 0\le\tau\le \taub$
being the minimizing ${\Cal L}$-geodesic joining $x_0$ to $y$.
Putting (A9)--(A11) together, one obtains a new proof of the
result of Cheeger-Yau, which asserts that if $M$ has nonnegative
Ricci curvature, then
$$
\heatau \left(\frac{e^{-\ell(y,
\tau)}}{(4\pi\tau)^{\frac{n}{2}}}\right)\le 0. \tag A12
$$
Namely $\frac{e^{-\ell(y, \tau)}}{(4\pi\tau)^{\frac{n}{2}}}$ is a
sub-solution of the heat equation. In particular,
$$
\frac{d}{d\tau} \int_M \frac{e^{-\ell(y,
\tau)}}{(4\pi\tau)^{\frac{n}{2}}}\, d\mu \le 0.
$$
Using the above geometric consideration, one can think the above
result of Cheeger-Yau as a parabolic volume comparison with the
respect to the positive measure $\frac{e^{-\ell(y,
\tau)}}{(4\pi\tau)^{\frac{n}{2}}}\, d\mu$. Recall that the
well-known Bishop-Gromov volume comparison states that if $M$ has
nonnegative Ricci curvature
$$
\frac{d}{d r} \left(\frac{1}{r^{n}}\int_{S_{x_0}(r)}dA\right)\le 0
$$
where $S_{x_0}(r)$ denotes the boundary of the geodesic ball
centered at $x_0$ with radius $r$, $dA$ is the induced area
measure. The by-now standard relative volume comparison can be
formulated in the similar way as above. Let $A$ be a measurable
subset of $S^{n-1}\subset T_{x_0}M$ one can define $C_A(r)$ to be
the collection of vectors $rv$ with $v\in A$ such that the
geodesic $\exp_{x_0}(sv)$ is minimizing for $s\le r$, where
$\exp_{x_0}(\cdot)$ is the (classical) exponential map. Then the
classical relative volume comparison theorem (cf. \cite{Gr})
asserts that if $M$ has nonnegative Ricci curvature
$$
\frac{d}{d r}\left(\frac{1}{r^n}\int_{\exp(C_A(r))}dA\right)\le 0.
$$
The following is a parabolic version of such relative volume
comparison theorem parallel to Perelman's work on Ricci flow
geometry.

\proclaim{Theorem A1} Assume that $M$ has nonnegative Ricci
curvature. Then
$$\frac{d}{d\tau}\int_{{\Cal
L}\exp_{C_A(\tau)}(\tau)}\frac{e^{-\ell(y,
\tau)}}{(4\pi\tau)^{\frac{n}{2}}}\, d\mu \le 0. \tag A13
$$
\endproclaim
\demo{Proof} The proof follows the similar argument as in
\cite{P}. Using the notation in the above discussion, first we
observe that
$$
\int_{{\Cal L}\exp_{C_A(\tau)}(\tau)}\frac{e^{-\ell(y,
\tau)}}{(4\pi\tau)^{\frac{n}{2}}}\,
d\mu=\int_{C_A(\tau)}\frac{e^{-\ell(y,
\tau)}}{(4\pi\tau)^{\frac{n}{2}}}J(\tau)\, d\mu_0,
$$
where $J(\tau)$ is the Jacobian of ${\Cal L}\exp_{(\cdot)}(\tau)$.
Since $C_A(\tau)$ is decreasing in $\tau$, it suffices to show
that
$$
\frac{d}{d \tau}\left(\frac{e^{-\ell(y,
\tau)}}{(4\pi\tau)^{\frac{n}{2}}}J(\tau)\right)\le 0. \tag A14
$$
This follows from (A9)--(A10), the fact that
$$
\frac{d}{d\tau} e^{-\ell}=\left(-\frac{\ell}{\taub}+\langle\nabla
\ell, X\rangle(\taub)\right)e^{-\ell}=0
$$
(since $\nabla \ell (\taub)=X$),  and the  claim that
$$
\frac{d}{d \tau} \log J(\taub) \le
\frac{n}{2\taub}-\frac{1}{\taub}\int_0^{\taub}\tau^{\frac{3}{2}}\Ric(X,
X)\, d\tau. \tag A15
$$
The estimate (A15) follows exactly the same proof as in \cite{P},
via the second variation formula (A6) and its consequence on the
Hessian comparison:
$$
\text{Hess}(L)(\bar{Y}, \bar{Y})\le
\int_0^{\taub}2\sqrt{\tau}\left(|\nabla_X Y|^2-R(X, Y, X,
Y)\right)\, d\tau
$$
where $Y(\tau)$ is a vector field along the minimizing ${\Cal
L}$-geodesic $\gamma(\tau)$ joining $x_0$ to $y$ with
$\gamma(\taub)=y$, satisfying   $Y(\taub)=\bar{Y}$ and
$$
\nabla_X Y-\frac{1}{2\tau}Y=0.
$$

\enddemo
Better comparison between Theorem A1 and the classical relative
volume comparison can be seen by comparing $\tM$ with $M$,
$M\times\{a\}$ with $S_{x_0}(a)$. Notice that in \cite{P}, one
does need such a result (not just the global version on the
monotonicity of $\tilde V_{g(\tau)}(x_0, \tau)$) to prove the
non-collapsing result on finite time solution to Ricci flow. This
and the potential application to the study of Riemannian geometry
of our linear version, justify the spelling out of this  result in
spite of its simplicity.

 \proclaim{Remark} One can formulate the
similar parabolic version of the volume comparison result for the
case $\Ric(M)\ge -(n-1)K$. We leave that to the interested
readers.
\endproclaim

In \cite{N3}, we showed that the lower bound on the entropy
$\inf_{0\le \tau\le T}\mu(\tau)$ (please see \cite{N3} for
definition) implies the non-collapsing of volume of ball of radius
$r$ for $r^2\le T$. This corresponds (but is considerably easier
than) the $\kappa$-non-collapsing result of Perelman on the
solutions to Ricci flow, via the monotonicity of entropy
functional. Using Theorem A1, one can have the following
consequence, which is just the linear version of the second proof
of Perelman on his celebrated $\kappa$-non-collapsing result.

\proclaim{Lemma A1} Let $M$ be a complete Riemannian manifold with
nonnegative Ricci curvature. Assume that $B(x_0, r)$ is
$\kappa$-collapsed in the sense that
$$
V(x_0, r)\le \kappa r^n. \tag A16
$$
Then there exists $C=C(n)>0$
$$
\tilde V(\kappa^{\frac{1}{n}}r^2)\le C\left(
\sqrt{\kappa}+\exp\left(-\frac{1}{8
\kappa^{\frac{1}{n}}}\right)\right). \tag A17
$$
Namely, the smallness of the relative volume $\frac{V(x_0,
r)}{r^n}$ is equivalent to  the smallness of the `reduced volume'
$\tilde V(r^2)$.
\endproclaim
Notice that the lemma  can also be proved by direct computations
and the Bishop volume comparison theorem, without using Theorem
A1.

\enddocument